\documentclass{aims}
\usepackage{amsmath}
  \usepackage{paralist}
  \usepackage{graphics} 
  \usepackage{epsfig} 
\usepackage{graphicx}  \usepackage{epstopdf}
 \usepackage[colorlinks=true]{hyperref}
\hypersetup{urlcolor=blue, citecolor=red}

\usepackage{lettrine}
\usepackage{tikz-cd}
\usepackage{amsthm,amsmath,amssymb}
\usepackage{graphicx}
\usepackage{subfigure}
\usepackage{changes}
\usepackage{soul}
\usepackage{tikz}
\usepackage{float}
\usepackage{tikz-cd}
\usetikzlibrary{arrows}
\usepackage{scalefnt}
\usepackage{makecell}
\usepackage[nomessages]{fp}
\usepackage{tikz-qtree}
\usetikzlibrary{shapes,positioning,chains,fit,calc}
\usepackage{rotating}
\usepackage{cancel}
\usepackage{algorithm,algorithmic}
\usepackage{multirow}
\usepackage{relsize}
\usepackage{pagecolor,lipsum}%

\usepackage{xspace}
\newcommand{\soc}{\textsf{SOC}\xspace}
\usepackage[section]{placeins}
\makeatletter
\AtBeginDocument{%
	\expandafter\renewcommand\expandafter\subsection\expandafter{%
		\expandafter\@fb@secFB\subsection
	}%
}
\makeatother

\definechangesauthor[name={Ilya}, color=red]{is}
\definechangesauthor[name={Hayato}, color=blue]{hm}

\AtBeginDocument{
	\label{CorrectFirstPageLabel}
	
}

\theoremstyle{definition}

\title[Optimal Placement of Wireless Charging Lanes] 
      {Optimal Placement of Wireless Charging Lanes in Road Networks}

\author[Ushijima-Mwesigwa,  Khan, Chowdhury, and Safro ]{}


 \email{hushiji@clemson.edu}
 \email{mdzadik@clemson.edu}
 \email{mac@clemson.edu}
 \email{isafro@clemson.edu}

\thanks{$^*$ Corresponding author: hushiji@clemson.edu}

\begin{document}
\maketitle

\centerline{\scshape Hayato Ushijima-Mwesigwa$^*$}
\medskip
{\footnotesize
 \centerline{ School of Computing}
   \centerline{Clemson University}
   \centerline{ Clemson, SC 29634, USA}
} 

\medskip

\centerline{\scshape MD Zadid Khan and Mashrur A. Chowdhury}
\medskip
{\footnotesize
 \centerline{ 	 Department of Civil Engineering}
   \centerline{Clemson University}
   \centerline{Clemson, SC 29634, USA}
}
\medskip

\centerline{\scshape Ilya Safro}
\medskip
{\footnotesize
 \centerline{ School of Computing}
	\centerline{Clemson University}
	\centerline{Clemson, SC 29634, USA}
}
\bigskip


\begin{abstract}
The emergence of electric vehicle wireless charging technology, where a whole lane can be turned into a charging infrastructure, leads to new challenges in the design and analysis of road networks. From a network perspective, a major challenge is determining the most important nodes with respect to the placement of the wireless charging lanes. In other words, given a limited budget, cities could face the decision problem of where to place these wireless charging lanes. With a heavy price tag, a placement without a careful study can lead to inefficient use of limited resources. In this work, the placement of wireless charging lanes is modeled as an integer programming problem. The basic formulation is used as a building block for different realistic scenarios. We carry out experiments using real geospatial data and compare our results to different network-based heuristics.\\
\noindent Reproducibility: all datasets, algorithm implementations and mathematical programming formulation presented in this work are available at\\
\noindent \url{https://github.com/hmwesigwa/smartcities.git}
\end{abstract}
\section{Introduction}\label{intro}
The transportation sector is the largest consumer in fossil fuel worldwide. As cities move towards reducing their carbon footprint, electric vehicles (EV) offer the potential to reduce both petroleum imports and greenhouse gas emissions. However, the batteries of these vehicles have a limited travel distance per charge. Moreover, the batteries require significantly more time to recharge compared to refueling a conventional gasoline vehicle. An increase in the size of the battery would proportionally increase the driving range. However, since the battery is the single most expensive unit in an EV, increasing its size would greatly increase the price. As a result leading to a major obstacle in EV widespread adaptation, \emph{range anxiety}, the persistent worry about not having enough battery power to complete a trip. 

Given the limitations of on-board energy storage, concepts such as battery swapping \cite{pan2010locating} have been proposed as possible approaches to mitigate these limitations. In the case of battery swapping, the battery is exchanged at a location that stores the equivalent replacement battery. This concept leads to issues such as battery ownership in addition to significant swapping infrastructure costs.
Another approach to increase the battery range of the EV is to enable power exchange between the vehicle and the grid while the vehicle is in motion. This method is sometimes referred to as \emph{dynamic charging} \cite{vilathgamuwa2015wireless,li2015optimizing} or \emph{charging-while-driving} \cite{chen2016optimal}. In this approach, the roads can be electrified and turned into charging infrastructure \cite{he2013integrated}. Dynamic charging is shown in \cite{ko2013optimal} to significanly reduce the high initial cost of EV by allowing the battery size to be downsized. This method could be used to complement other concepts such as battery swapping to reduce driver range anxiety.

There have been many studies on the design, application and future prospects of wireless power transfer for electric vehicles (see e.g., \cite{qiu2013overview,bi2016review,li2015wireless,lukic2013cutting, cirimele2014wireless,fuller2016wireless,vilathgamuwa2015wireless,ning2013compact}). Some energy companies are teaming up with automobile companies to incorporate wireless charging capabilities in EVs. Examples of such partnerships include Tesla-Plugless and Mercedez-Qualcomm. Universities, research laboratories and companies have invested in research for developing efficient wireless charging systems for electric vehicles and testing them in a dynamic charging scheme. Notable institutions include Auckland University \cite{cho2011magnetic}, HaloIPT (Qualcomm) \cite{lee2010line}, Oak Ridge National laboratory (ORNL) \cite{kim2012design}, MIT (WiTricity) and Delphi \cite{jung2012high}. However, there is still a long way to go for a full commercial implementation, since it requires significant changes to be made in the current transportation infrastructure.

	A few studies focus on the financial aspect of the implementation of a dynamic charging system.  A smart charge scheduling model is presented in \cite{li2015optimizing} that maximizes the net profit to each EV participant while simultaneously satisfying energy demands for their trips. An analysis of the costs associated with the implementation of a dynamic wireless power transfer infrastructure and a business model for the development of a new EV infrastructure are presented in \cite{gill2014infrastructure}. Integrated pricing of electricity in a power network and usage of electrified roads in order to maximize the social welfare is explored in \cite{he2013integrated}.
    
In regards to the planning infrastructure, a number of studies have focused on the implications of dynamic charging to the overall transportation network. An analysis on the effectiveness of placing wireless charging units at traffic intersections in order to take advantage of the frequent stops at these locations is taken in \cite{mohrehkesh2011toward} . Methods on how to effectively distribute power to the different charging coils along a wireless charging lane in a \emph{vehicle-to-infrastructure} (V2I) communication system have also be demonstrated \cite{sarker2016efficient}. The authors in 
\cite{johnson2013utilizing} carry out simulations over a traffic network to show how connected vehicle technology, such as \emph{vehicle-to-vehicle} (V2V) or V2I communications can be utilized in order effectively facilitate the EV charging process at fast-charging stations.
Routing algorithms that take dynamic charging into account have also be developed. An ant colony optimization based multi-objective routing algorithm that utilizes V2V and V2I communications systems to determine the best route considering the current battery charge is developed in \cite{li2016connectivity}.

Given the effectiveness and advances in dynamic charging technology, cities face the challenge of budgeting and deciding on what locations to install these \emph{wireless charging lanes} (WCL) within a transportation network. In this article, we seek to optimize the installation locations of WCLs.

\subsection{Related Work}
Owing to advances in technology, there have been recent studies related to the optimal placement of wireless charging lanes. The basic difference in these studies arise in the objective function and/or the type of routes, between the origin and destination, that are considered.

In a recent study \cite{chen2016optimal}, the optimal placement of wireless charging lanes when the charging infrastructure is considered to affect the EV driver's route choice is developed. They developed a mathematical model with an objective to minimize the total system travel times which they defined as the total social cost.
There have also been studies devoted optimal locations of refueling or recharging stations of EVs when the EV driver route choice is not fixed (see, e.g, \cite{he2013optimal,kang2014strategic,jung2014stochastic,he2015deploying}). 

One prominent study where the charging infrastructure does not affect the route choice is presented in \cite{ko2013optimal}. In this study, they focus on a single route and seek the optimal system design of the \emph{online electric vehicle} (OLEV) that utilizes wireless charging technology. They apply a \emph{particle swarm optimization} (PSO) method to find a minimum cost solution considering the battery size, total number of WCLs (power transmitters) and their optimal placement as decision variables. The model is calibrated to the actual OLEV system and the algorithm generates reliable solutions. However, the formulation contains a non-linear objective function making it computationally challenging for multi-route networks. Moreover, speed variation is not considered in this model, which is typical in a normal traffic environment. The OLEV and its wireless charging units were developed in Korea Advanced Institute of Science and Technology (KAIST) \cite{jang2012optimal}. At Expo 2012, an OLEV bus system was demonstrated, which was able to transfer 100KW ($5\times 20$KW pick-up coils) through 20 cm air gap with an average efficiency of 75\%. The battery package was successfully reduced to 1/5 of its size due to this implementation \cite{jang2012creating}. This study was recently extended \cite{mouhrim2016optimal} to take multiple routes into account in which they carried out experiments on example with five routes.

In the literature, studies on optimal locations of plug-in charging facilities are often related to the maximal covering location problem (MCLP), in which each node has a demand and the goal is to maximize the demand coverage by locating a fixed number of charging facilities. For a more comprehensive study on the MCLP, one can refer to the work in \cite{church1974maximal,farahani2012covering,daskin2008you,hale2003location}. The flow-capturing location problem (FCLP) \cite{hodgson1990flow} builds on the MCLP and defines which seeks to maximize the \emph{captured flow} between all origin-destination pairs. Flow along a path is defined as being captured if there exists at least one facility on the path. The definition of a flow being captured however does not carry over to the case of vehicle refueling, as a vehicle may need to refuel more than once to successfully complete the entire path. As a result, the flow-refueling location model (FRLM) is formulated in \cite{kuby2005flow}. Subsequently, extensions of the FRLM have been formulated (see, e.g., \cite{kuby2007location,upchurch2009model,kim2012deviation,huang2015optimal,kim2013network,wang2013locating}). In the FRLM and its extensions, the assumption that the vehicle is fully refueled at a facility does not carry over to the case of in-motion wireless charging as a EV may not be fully charged after passing over a wireless charging unit. An extention of FRLM  where the routes are not fixed is given in \cite{riemann2015optimal} where they apply their formulation to wireless charging facilities. Similarly to other FRLM extensions, they assume that an EV is fully charged once it passes over a link containing a wireless charging facility. The flow-based set covering model for fast-refueling stations such as battery exchange or hydrogen refueling stations is proposed in \cite{wang2009locating}. In particular, their approach does not assume that the fuel or charge after passing through refueling or recharging facility to be full. Their work was subsequently extended \cite{wang2010locating,wang2011locating} while keeping a similar objective to minimize the locating cost. 

A different approach is taken in \cite{dong2014charging} in order to optimize the locations of public charging facilities for EVs. They take into account the long charging times of these charging stations which increases the preference for a charging facility to be located at a user activity destination.  
\subsection{Contribution}
In this work, we  seek to address the optimal location problem of wireless charging lanes in road networks, given a limited budget. Our objective is to maximize the number of origin-destination routes that benefit, to a given threshold, from a deployment. This objective function is different from the one considered in \cite{chen2016optimal}.
Given that the battery charge may not significantly increase when an EV drives over a single wireless charging lane, the minimum budget to cover an entire network may be significantly higher than the available budget. In order to best utilize the available budget, we define a feasible path, as an origin-destination path,  in relation to the final battery charge an EV would have at the end of its trip along this path. We then seek to maximize the number of feasible paths over the network.  We formulate the WCL installation problem as an integer programming model that is built upon taking into account different realistic scenarios. We compare the computational results for the proposed model to  faster heuristics and demonstrate that our approach provides significantly better results for fixed budget models. Using a standard optimization solver with parallelization, we provide solutions for networks of different sizes including the Manhattan road network, whose size is significantly larger than the ones considered in previous studies.


	\section{Optimization Model Development}\label{optmodel}
	The purpose of developing a mathematical model of the WCL installation problem is to construct an optimization problem that maximizes the battery range per charge within a given budget and road network. This, in turn, will minimize the driver range anxiety within the road network. In this section, relevant definitions followed by modeling assumptions are presented.
	\subsection{Road Segment Graph}
    Consider a physical network of roads within a given location. A road segment is defined as the one-way portion of a road between two intersections. Let $G = (V,E)$ be a directed graph with node set $V$ such that $v \in V$ if and only if $v$ is a road segment. Two road segments $u$ and $v$ are connected with a directed edge $(u,v)$ if and only if the end point of road segment $u$ is adjacent to the start point of road segment $v$. We refer to this graph as a \emph{road segment graph}. This representation is adopted over the conventional network representation because the decision variables are based on road segments. Other advantages to this representation such as modeling of turn costs for a given route are for example given in \cite{bogaert2005line, Winter2002}. For a given road segment graph and budget constraint, the objective is to find a set of nodes that would minimize driver range anxiety within the network. 

\subsection{State of Charge of an EV}\label{soc_section}

State of charge ($\soc$) is the equivalent of a fuel gauge for the battery pack in a battery electric vehicle and hybrid electric vehicle. 
 The $\soc$ determination is a complex non-linear problem and there are various techniques to address it (see e.g., \cite{chang2013state,wang2016comparison,koirala2015comparison}). As discussed in the literature, the $\soc$ of an EV battery can be determined in real time using different methods, such as terminal voltage method, impedence method, coulomb counting method, neural network, support vector machines, and Kalman fitering. The input to the models are physical battery parameters, such as terminal voltage, impedence, and discharging current. However, the SOC related input to our optimization model is the change in $\soc$ of the EV battery to traverse a road segment rather than the absolute value of the real time $\soc$ of the EV battery. So, we formulate a function that approximates the change in $\soc$ of an EV to traverse a road segment using several assumptions, as mentioned in the following.
The units of $\soc$ are assumed to be percentage points (0\% = empty; 100\% = full). The change in $\soc$ is assumed to be proportional to the change in battery energy. This is a valid assumption for very small road segments that form a large real road network, which is the case in this analysis (range of 0.1 to 0.5 mile). 

We compute the change in $\soc$ of an EV as a function of the time $t$ spent traversing a road segment by
\begin{equation}
\Delta \soc_{t} = \frac{E_{\text{end}} - E_{\text{start}}}{E_{\text{cap}}},
\label{soc_change}
\end{equation}
where $E_{\text{start}}$ and $E_{\text{end}}$ is the energy of the battery (KWh) before and after traversing the road segment respectively and $E_{\text{cap}}$ is the battery energy capacity. We follow the  computation of the battery energy given by \cite{sarker2016efficient}. We, however, assume that the velocity of an EV is constant while traversing the road segment. This gives us 
\begin{equation}
E_{\text{end}} - E_{\text{start}} =   (P_{2t}\cdot \eta)t - P_{1t}t ,
\label{energy}
\end{equation}
where $P_{1t}$ is the power consumption (KW) needed to traverse the given road segment in time $t$. $P_{2t}$ is the power delivered to the EV in case a WCL is installed on the road segment, otherwise  $P_{2t}$ is zero.  In order to take into account the inefficiency of charging due to factors such as misalignment between the primary (WCL) and secondary (on EV) charging coils and air gap, an inefficiency constant $\eta$ is assumed.
 
The power consumption $P_{1t}$ varies from EV to EV. In this work, we take an average power consumption calculated by taking the average mpge (miles per gallon equivalent) and battery energy capacity rating from a selected number of EV. We took the average of over 50 EVs manufactured in 2015 or later. For each EV, its fuel economy data was obtained from \cite{fuelecon}. The values of $P_{2t}$ and $\eta$, the power rating of the WCL, and the efficiency factor, we average the values from \cite{bi2016review}, Table 2, where the authors make a comparison of prototype dynamic wireless charging units for electric vehicles.

\subsubsection{Reliability Of \soc}
In reality, an accurate estimation of the change in \soc depends on different factors that are not considered in this work such as acceleration/deceleration,  and elevation of road segment  \cite{fontana2013optimal}. However, assumptions on the estimation and variation of \soc must be made for the optimization problem formulation. Similar studies such as \cite{he2013integrated,chen2016optimal} use a linear model based only on distance to the estimate change of \soc, while \cite{riemann2015optimal} simply assumes that an EV will be fully charged if it passes over a charging lane. In this work, we estimate the variation of \soc using the nonlinear function given in the previous section.

The estimation of \soc of an EV can be a difficult task compared to Internal Combustion Engines (ICE) where the fuel level in the tank is monitored by sensors~\cite{watrin2012review}.
Traditional \soc estimation methods for EVs are open loop methods such as the ampere-hour integral method \cite{zhang2015integrated}. In the ampere-hour integral method, the current is integrated over a certain time interval to determine the \soc value. Due to the need for sampling of current values for integral calculation, a sampling error is introduced, and the sampling errors are accumulated over time. Furthermore, in order to determine the final \soc value from the integral, the initial \soc value is required, which is also difficult to determine \cite{sun2016systematic}. Closed loop systems have a feedback loop that solves the aforementioned issues, so advanced filtering techniques with feedback loops are coupled with an EV battery model to form a more accurate SOC estimation model. The output from the EV battery model is a voltage that is compared with the actual / measured output voltage of the battery. The difference between the two voltages is fed back to the system, forming a closed loop system. The accuracy of the \soc estimation model is primarily determined by the accuracy of the EV battery model. There are two types of EV battery models; the equivalent circuit models and the electro-chemistry models. In the equivalent circuit models, the battery is represented as a combination of a DC voltage source, resistors and capacitors \cite{greenleaf2015temperature}. The drawback of these models is the difficulty in identifying the structure and the parametric values of the circuit that accurately represents the battery. The electro-chemistry based models mathematically formulate the chemical reactions in the battery, so the model structure and its parameters are easier to determine, but the high computational complexity of the models make it unsuitable for real time SOC estimation \cite{stetzel2015electrochemical}. Moreover, other factors such as temperature, depth of discharge, cycling, recharging voltage and maintenance have high impact on the chemical reactions in the batteries, so all these factors are considered in the electro-chemistry based models. From the above discussion, it can be concluded that the current \soc estimation methods still have some inaccuracies and uncertainties associated with it, and research is ongoing to develop more accurate \soc estimation models.

\subsection{Modeling Assumptions}\label{assumptions}

For a road segment graph $G$, we assume that each node has attributes such as average speed and distance that are used to compute the average traversal time of the road segment. Since nodes represent road segments, an edge represents part of an intersection, thus, the weight of an edge does not have a typical general purpose weighting scheme associated to it (such as a length). Since the proposed model is developed to optimize WCL placement for a given set of routes within a road network, in computational experiments, routes are chosen based on travel time. 
As a result, each edge $(u,v)$ in $G$ is assigned a weight equal to the average traversal time of road segment $u$.

In computational experiments, we assume that $\soc$ of any EV whose journey starts at the beginning of a given road segment is fixed. (However, this assumption can easily dropped with a little modification of the model if real information about initial SOC is available.) For example, we may assume that if a journey starts at a residential area, then any EV at this starting location will be fully charged  or follows a charge determined by a given probability distribution which would not significantly change the construction of our model. For example, in real applications, one could choose the average $\soc$ of EV's that start at that given location. In our empirical studies, for simplicity, we first assume that all EVs start fully charged. Results for studies where the initial \soc is chosen uniformly at random are also given. We
assume that $\soc$ takes on real values such that $0 \leq \soc \leq 1$ at any instance where $\soc =1$ implies that the battery is fully charged and $\soc=0$ implies that the battery is empty.

Our model is based on optimizing the WCL placement with respect to a set of routes. In the computational experiments, we assume that longer routes will have higher priority for the WCL placement and thus we randomly choose longer routes to be considered as input to the model. We assume that it is not necessary to consider all routes. This is because some routes could significantly affect the outcome of the model, however in reality a large number of them could be routes with very low priority for WCL placement. For example, these can be routes that users usually have sufficient charge to complete, such as a route from a users home to a local grocery store. Therefore, we do not include short distance routes in our experiments; on the other hand, we do not have information on what medium and long distance routes should be considered. Moreover, we cannot know it now with a good level of certainty at least because, in general, EVs and autonomous vehicles (that are expected to change the ways we use vehicles and roads) are still not dominating the market (not to mention EVs that use WCL). Thus, the choice of longer routes chosen at random is a way to demonstrate the use of our model.  

A route is \emph{infeasible} within a network if any EV that starts its journey at the beginning of this route (starts fully charged in our empirical studies), will end with a final $\soc \leq \alpha$, where $0 \leq \alpha \leq 1$. The constant $\alpha$ is a global parameter of our model called a global $\soc$ threshold. The value of $\alpha$ could be chosen in relation to the minimum $\soc$ an EV driver is comfortable driving with \cite{franke2012experiencing}. Introducing different types of EV and more than one type of $\alpha$ would not significantly change the construction of the model. 

Given the total length of all road segments in the network, $T$, we define the budget, $0\leq \beta \leq 1$, as a fraction of $T$ for which funds available for WCL installation. 
For example, if $\beta=0.5$, the city planners have enough funds to install WCL's across half the length of the entire road network. In this research, the model and its variations are used to answer the following problems that the city planners are interested in.

\begin{enumerate}
	\item For a given  $\alpha$, determine the \emph{minimum budget}, $\beta$, together with the corresponding locations, needed such that the number of infeasible routes is zero.
	\item For a given $\alpha$ and $\beta$, determine the optimal installation locations to \emph{minimize the number of infeasible routes}.
\end{enumerate}

We assume that minimizing the number of infeasible routes would reduce the driver range anxiety within the network.

	\subsection{Single Route Model Formulation}
Let $Routes$ be a set of routes in the road segment graph $G$.
For each route, $r \in Routes$, assume that each EV whose journey is identical to this route has a fixed initial $\soc$, and a variable final $\soc$, termed $i\soc_r$, and $f\soc_r$, respectively, depending on whether or not WCL's were installed on any of the road segments along the route. The goal of the optimization model is to guarantee that either $f\soc_r \geq \alpha$ where $\alpha $ is a global threshold or $f\soc_r$ is as close as possible to $\alpha $ for a given budget. Given that realistic road segment graphs have a large number of nodes, taking all routes into account may overwhelm the computational resources, thus, the model is designed to give the best solution for any number of routes considered. 

The proposed model is first described for a single route and then generalizing it to multiple routes. For simplicity, we will assume that the initial $\soc$, $i\soc_r  =1$, for each route $r$, i.e., all EV's  start their journey fully charged. This assumption can easily be adjusted with no significant changes to the model.  For ease of exposition, in this section we will also assume a  simplistic $\soc$ function to estimate \soc levels, that is, $\soc$ of an EV traversing a given road segment increases by one discrete \soc level if a WCL is installed, otherwise it decreases by one discrete \soc level. In the experiments based on real data, we use actual energy consumption and energy charged based on the change of \soc described by equations (\ref{soc_change}) and (\ref{energy}).

For a single route $r \in Routes$, with $i\soc_r =1$, consider the problem of determining the optimal road segments to install WCL's in order to maximize $f\soc_r$ within a limited budget constraint. Define a \soc-state graph, $\mathcal{G}_r=(\mathcal{V}_r, \mathcal{E}_r)$, for route $r$, as an acyclic directed graph whose vertex set, $\mathcal{V}_r$, describes the varying $\soc$ an EV on a road segment would have depending on whether or not the previously visited road segment had a WCL installed. More precisely, let $r = (u_1, u_2, \dots, u_m)$, for $u_i \in V$ with $i = 1, \dots, m$ and $m >0$. Let $nLayers \in \mathbb{N}$ represent the number of discrete values that the $\soc$ can take. 
For each $u_i \in r$, let $\mu_{i, j} \in \mathcal{V}_r$ for $j=1, \dots, nLayers$ representing the $nLayers$ discrete values that the $\soc$ can take at road segment $u_i$.
Let each node $\mu_{i,j}$ have out-degree at most 2, representing the two different scenarios of whether or not a WCL is installed at road segment $u_i$. An edge in the edge set, $\mathcal{E}_r$, will be represented by a triple $(r,\mu_{i_1,j_1},\mu_{i_2,j_2})$ where $\mu_{i_1,j_1}, \mu_{i_2, j_2} \in \mathcal{V}_r$. The edge $(r, \mu_{i, j_1}, \mu_{i+1, j_2})$ is assigned weight 1 to represent the scenario if a WCL is installed at $u_i$ and 0, otherwise. An extra node is added accordingly to capture the output from the final road segment $u_m$, we can think of this as adding an artificial road segment $u_{m+1}$. Two dummy nodes, $s$ and $t$, are also added to the \soc-state graph $\mathcal{G}_r$ to represent the initial and final $\soc$ respectively. There is one edge of weight 0 between $s$ and $\mu_{1,j^*}$ where the node $\mu_{1,j^*}$ represents the initial $\soc$ of an EV on this route. Each node $\mu_{m+1,j}$ for each $j$ is connected to node $t$ with weight 0. 

\begin{figure}[ht]
\centering
\includegraphics[clip, trim=4.5cm 16.2cm 3cm 4cm,width=.7\textwidth]{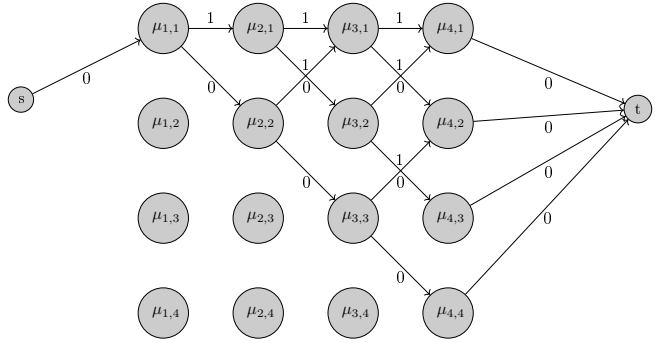}
\caption{ Example of $\mathcal{G}_r$ with $r = (u_1, u_2, u_3)$ and $nLayers = 4$. $u_4$ is an artificial road segment added to capture the final $\soc$ from $u_3$. The nodes in the set $\mathcal{B}_r =\{\mu_{i,j} | i=4 \text{ or } j =4\}$ are referred to as the boundary nodes. The out going edges of each node $\mu_{i,j}$ are determined by an $\soc$ function.  Each node represents a discretized $\soc$ value.}
\label{soc_graph}
\end{figure}
Consider a path $p$ from $s$ to $t$, namely, $p = (s, \mu_{1,j_1}, \mu_{2, j_2}, \dots, \mu_{m, j_m}, t)$, then each node in $p$ represents the \soc of an EV along the route. We use this as the basis of our model. Any feasible $s$-$t$ path will correspond to an arrival at a destination with an $\soc$ above a given threshold. A minimum cost path in this network would represent the minimal number of WCL installations in order to arrive at the destination.  Figures \ref{soc_graph} shows an example for a $\soc$-state graph constructed from a single route with three road segments $u_1, u_2$ and $u_3$ with four discretized \soc levels taken into account. The nodes $\mu_{4,j}$ for $j =1,\dots, 4$ are added to capture the output \soc level from road segment $u_3$. The nodes $\mu_{i,1}$ and $\mu_{i,4}$, for $i = 1, \dots, 4$, represent the maximum and minimum \soc levels respectively, for the road segments $u_i$.  The $s$-$t$ path $p=(s, \mu_{1,1}, \mu_{2,2}, \mu_{3,1}, \mu_{4,2}, t)$ would, for example, give an \soc level if a WCL was installed on road segment $u_2$.  
If $\upsilon, \nu \in \mathcal{V}_r$, with $(r, \upsilon, \nu) \in \mathcal{E}_r$ and weight $w_{r,\upsilon,\nu}$, then the minimum cost path can be formulated as follows: 
	\begin{equation}
	\begin{array}{ll@{}ll}
	\text{minimize}  & \displaystyle\sum\limits_{(r, \upsilon, \nu) \in \mathcal{E}_r}w_{r,\upsilon,\nu}x_{r,\upsilon, \nu} &\\
	\text{subject to}& \displaystyle\sum\limits_{\nu \in \mathcal{V}_r} x_{r,\upsilon,\nu} - \displaystyle\sum\limits_{\nu \in \mathcal{V}_r} x_{r,\nu,\upsilon} = \begin{cases}
1, &\text{if } \upsilon  = s;\\
-1,  &\text{if }\upsilon =t;\\
0, & \text{otherwise}
	\end{cases} \quad \forall \upsilon \in \mathcal{V}_r \\
	&   x_{r,\upsilon,\nu} \in \{0,1\}
	\end{array}
    \label{mincost}
	\end{equation}
	where $\displaystyle\sum\limits_{\nu \in \mathcal{V}_r} x_{r,\upsilon,\nu} - \displaystyle\sum\limits_{\nu \in \mathcal{V}_r} x_{r,\nu,\upsilon} = 0$ ensures that we have a path i.e., number of incoming edges is equal to number of out going edges.

\paragraph{Decision Variable for Installation}
	Let $R_k$ be the decision variable for installation of a WCL at road segment $u_k$. Then, for a single route, we have $R_k = w_{r,\upsilon^*,\nu^*} \in \{0,1\}$, where $(r, \upsilon^*, \nu^*) \in \mathcal{E}_r$ is an edge belonging to the minimum cost path of the optimal solution of (\ref{mincost}), i.e., $x_{r, \upsilon^*, \nu^*} =1$, with $\upsilon^* = \mu_{k,i}$ for some  $i$, and $\nu^* = \mu_{k+1, j}$, for some $j$. 
    Under the constraints for a single route, an optimal solution to the minimum cost path from $s$ to $t$ would be a solution for the minimum number of WCL's that need to be installed, in order for EV to arrive at the destination with its final $\soc$ greater than a specified threshold. 	
	
	\subsection{General Model Description}\label{description}
    
   Consider the case with multiple routes. For each route, $r\in Routes$, an \soc-state graph, $\mathcal{G}_r = (\mathcal{V}_r, \mathcal{E}_r)$ is formulated. Notice that since a road segment can belong to multiple routes, the set $\mathcal{V}_{r_1}\cap \mathcal{V}_{r_2}$ is not necessarily empty for two distinct routes, $r_1$ and $r_2$, however, $\mathcal{E}_{r_1}\cap \mathcal{E}_{r_2} = \emptyset$. 
   

	\noindent \textbf{Decision Variables:}\\
For a given route, $r$, and $\soc$-state graph $\mathcal{G}_r = (\mathcal{V}_r,\mathcal{E}_r)$ where the edges $E_r$ are defined according to an $\soc$ function, for example, the \soc function given by equation (\ref{soc_change}).
The weight of an edge $(r, \upsilon, \nu) \in \mathcal{E}_r$ for $\upsilon = \mu_{k,i} \in \mathcal{V}_r$, for some $i$, is given by
\begin{equation*}
	w_{r,\upsilon,\nu} = \begin{cases}
		1, & \text{if WCL is installed in respective road segment for $u_k$}  \\
		0, & \text{otherwise }
	\end{cases}
\end{equation*}
%
Then the decision variables of the model are given by
		\begin{align*}
		R_{k} &= \begin{cases}
		1, & \text{if at least one route requires a WCL installation at $u_k$}  \\
		0, & \text{otherwise }
		\end{cases}\\
	&&	\text{for } k = 1,\dots, nRoadSegs\\
	x_{r,\upsilon, \nu} &= \begin{cases}
	1, & \text{if edge $(r, \upsilon, \nu)$ is in an $s$-$t$ path in $G_r$}  \\
	0, & \text{otherwise }
	\end{cases}\\
		&&	\text{for } r = 1,\dots, nRoutes\\
		\end{align*}
	For the decision variable $R_k$ on the installation of a WCL at road segment $u_k$, we install a WCL if at least one route requires an installation within the different $s$-$t$ paths for each route. For road segment $u_k$, and for any set of feasible $s$-$t$ paths, let $p(u_k)$ be the number of routes that require a WCL installation at road segment $u_k$, then $p(u_k)$ is given by
		\begin{equation*}
		p({u_k}) =
		\displaystyle\sum\limits_{{r} = 1}^{nRoutes}
        \displaystyle\sum\limits_{\substack{(r, \upsilon,\nu)\in \mathcal{E}_r \\ \upsilon=\mu_{k,i}\\ i \in \mathbb{N}}}  w_{{r},\upsilon,\nu}\cdot x_{{r},\upsilon,\nu}
		\end{equation*}
		Then,
			\begin{equation}
			R_{{k}} = \begin{cases}
			1, & \text{if }  p(u_k)\geq 1 \\
			0, & \text{otherwise }
			\end{cases}
            \label{install_logic}
			\end{equation}
models the installation decision. 

\noindent \textbf{Objective function:}
    For the problem of minimizing the budget, the objective function is simply given by minimizing
    \begin{equation}
    \sum_{{k} =1 }^{nRoads} c_{{k}}\cdot R_{{k}}.
    \end{equation}
	where $c_k$ is the cost of installing a WCL at road segment $u_k$.
    
    For the problem of minimizing the number of infeasible routes for any fixed budget, we modify the $\soc$-state graph such that there exists an $s$-$t$ path for any budget. We achieve this by adding an edge of weight 0 between the nodes $\mu_{i, nLayers}$ to $t$ for all $i = 1, \dots, m + 1$, in each route, where $m$ is the number of road segments in the route. Define the \emph{boundary nodes} of the $\soc$-state graph with respect to route $r$, as the set of all the nodes adjacent to node $t$. Let $\mathcal{B}_r$ be the boundary nodes with respect to route $r$. Assign each node in $\mu_{i,j} \in \mathcal{B}_r$ weights according to the function:

	 \begin{equation*}
	\textbf{w}(\mu_{i,j}) = \begin{cases}
	 1, & \text{if } \mathbf{s}(\mu_{i,j}) \geq \alpha \\
	 0, & \text{otherwise}
	 \end{cases}
	 \end{equation*}
where $\mathbf{s}(\mu_{i,j})$ is the discretized $\soc$ value that node $\mu_{i,j}$ represents. 
In the weighting scheme above, there is no distinctions between two infeasible routes. However, a route in which an EV completes, say, 90\% of the trip would be preferable to one in which an EV completes, say, 10\% of the trip. This preference is taken into account and the weight the boundary nodes can also be given by the function
	 \begin{equation}
	 \textbf{w}(\mu_{i,j}) = \begin{cases}
	 1, & \text{if } \mathbf{s}(\mu_{i,j}) \geq \alpha \\
	  \frac{ d(r, u_i) - |r|}{|r|}, & \text{otherwise},
	 \end{cases}
     \label{boundary_weight}
	 \end{equation}
where $|r|$ is the distance of the route and $d(r,u_i)$ is the distance between the origin and the end of road segment $u_i$ for route $r$. The term $\frac{ d(r, u_i)-|r|}{|r|}$ is a penalty depending on how close an EV comes to completing a given route. 

In order to take the inaccuracies of determining the variation of \soc into account, we introduce a tolerance parameter $\epsilon_{\text{tol}}$  to the model such that a route with final \soc in the range $(\alpha - \epsilon_{\text{tol}}, \alpha + \epsilon_{\text{tol}})$, is not necessarily labeled as feasible/infeasible. An objective function is formed where such routes have larger contributions to the objective value compared to routes strictly less than $\alpha - \epsilon_{\text{tol}}$, however, contribute less to the objective value compared to routes with final \soc greater than $\alpha + \epsilon_{\text{tol}}$. To take this into account, the weight of the boundary nodes can be given by
	 \begin{equation}
	 \textbf{w}(\mu_{i,j}) = \begin{cases}
	 1, & \text{if } \mathbf{s}(\mu_{i,j}) \geq \alpha + \epsilon_{\text{tol}} \\
     0, & \text{if } \mathbf{s}(\mu_{i,j}) \in (\alpha - \epsilon_{\text{tol}}, \alpha + \epsilon_{\text{tol}}) \\
	  \frac{ d(r, u_i) - |r|}{|r|}, & \text{otherwise}.
	 \end{cases}
     \label{boundary_weight2}
	 \end{equation}
 
   Since two routes are not necessarily equal, that is, planners may not care about routes with low demand, we incorporate a normalized parameter $\delta_r \in (0,1)$ for each route $r$. The value of $\delta_r$ represents the normalized travel demand for the specific route. 
The objective is then given by maximizing the expression
  		\begin{equation}
		\displaystyle\sum\limits_{{r}=1}^{nRoutes}  \displaystyle\sum\limits_{\nu = \mu_{i,j}\in \mathcal{B}_{{r}}} \delta_r \cdot \textbf{w}(\nu) \cdot x_{r,\nu,t}.
		\end{equation}
	\noindent	\textbf{Budget Constraint:}\\
	Cost of installation cannot exceed a budget $B$. Since this technology is not yet widely commercialized, we can discuss only the estimates of the budget for WCL installation. Currently, the price of installation per kilometer ranges between a quarter million to several millions dollars \cite{gill2014infrastructure}. For simplicity, in our model the cost of installation at a road segment is assumed to be proportional to the length of the road segment which is likely to be a real case. Thus, a budget would represent a fraction of the total length of all road segments.
	\begin{equation}
	\sum_{{k} =1 }^{nRoads} c_{{k}}\cdot R_{{k}} \leq B
	\end{equation}
	
	\noindent \textbf{$s$-$t$ path constraints for \soc-state graph:}
The constraints defining an $s$-$t$ path for all $\upsilon \in \mathcal{V}_r$
	\begin{equation*}
	\displaystyle\sum\limits_{\nu \in \mathcal{V}_r} x_{r,\upsilon,\nu} - \displaystyle\sum\limits_{\nu \in \mathcal{V}_r} x_{r,\nu,\upsilon} = \begin{cases}
	1, &\text{if } \upsilon  = s;\\
	-1,  &\text{if } \upsilon =t;\\
	0, & \text{otherwise}
	\end{cases}
	\end{equation*}
    Note: In the current formulation of the proposed model, constraints for route feasibility are not explicitly needed because the \soc-state graph constructed takes this into account. 
\subsubsection{Model}

The complete model formulation for minimizing the number of infeasible routes, for a fixed budget is given by:
		\begin{align*}
		\begin{array}{ll@{}ll}
		\text{maximize}  & 	\displaystyle\sum\limits_{{r}=1}^{nRoutes}  \displaystyle\sum\limits_{\nu = \mu_{i,j}\in \mathcal{B}_{{r}}}\delta_r\cdot \textbf{ w}(\nu)\cdot x_{r,\nu,t} \\
		\text{subject to}&
		\sum_{{k} =1 }^{nRoads} c_{{k}}\cdot R_{{k}} \leq B\\
		&\displaystyle\sum\limits_{\nu \in \mathcal{V}_r} x_{r, \upsilon,\nu} - \displaystyle\sum\limits_{\nu \in \mathcal{V}_r} x_{r, \nu,\upsilon} = \begin{cases}
		1, &\text{if } \upsilon =  s\\
		-1,  &\text{if }  \upsilon =t\\
		0, & \text{otherwise}
		\end{cases}& \quad  r = 1, \dots, {nRoutes}, \ \upsilon \in \mathcal{V}_r \\
		&  R_{k}\leq p( u_k)&  {k} = 1, \dots, {nRoadSegs}\\
		& \mathbf{M}\cdot R_{k}\geq p( u_k)& {k} = 1, \dots, {nRoadSegs}\\
		&R_{k}  \in \{0,1\} &  k= 1, \dots, {nRoadSegs}
		\end{array}
		\end{align*}
		where 
\begin{equation*}
		p({u_k}) =
		\displaystyle\sum\limits_{{r} = 1}^{nRoutes}
        \displaystyle\sum\limits_{\substack{(r, \upsilon,\nu)\in \mathcal{E}_r \\ \upsilon=\mu_{k,i}\\ i \in \mathbb{N}}}  w_{{r},\upsilon,\nu}\cdot x_{{r},\upsilon,\nu}
		\end{equation*}
and $\mathbf{M}$ is a large constant used to model the logic constraints given in equation (\ref{install_logic}).

Since we are interested in reducing the number of routes with a final $\soc$ less than $\alpha $, we can take the set of routes in the above model to be all the routes that have a final $\soc$ below the given threshold. We evaluate this computationally and compare it with several fast heuristics.

	\subsection{Heuristics}\label{heuristics}
	Integer programming is NP-hard in general and since the status of the above optimization model is unknown, we have little evidence to suggest that it can be solved efficiently. For large road networks, it may be desirable to use heuristics instead of forming the above integer program. In particular, since we know the structure of the network, one natural approach may be to apply concepts from network science to capture the features of the best candidates for a WCL installation. In this section, we outline different heuristics for deciding on the set of road segments. We then compare these structural based solutions to the optimization model solution, and demonstrate the superiority of proposed model.
	
	Different centrality indexes is one of the most studied concepts in network science \cite{Newman:2010:NI}. Among them, the most suitable to our application are betweenness and vertex closeness centralities. In \cite{freeman1978centrality}, a node closeness centrality is defined as the sum of the distances to all other nodes where the distance from one node to another is defined as the shortest path (fastest route) from one to another. Similar to interpretations from \cite{borgatti1995centrality}, one can interpret closeness as an index of the expected time until the arrival of something "flowing" within the network. Nodes with a low closeness index will have short distances from others, and will tend to receive flows sooner. In the context of traffic flowing within a network, one can think of the nodes with low closeness scores as being well-positioned or most used, thus ideal candidates to install WCL.
	
\begin{figure}
    \centering
	\includegraphics[width=0.4\linewidth]{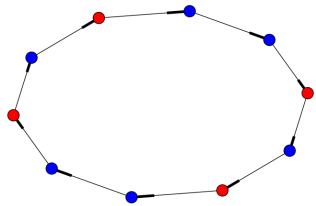}
    \caption{Optimal solution with a four unit installation budget. The thick ends of the edges are used to indicate the direction of the edge. 
  Taking $\alpha=0$ without any installation, there are 70 number of infeasible routes. An optimal installation of 5 WCLs would ensure zero infeasible routes. With an optimal installation of 4 WCLs, the nodes colored in red, there would have 12 infeasible routes.}
    \label{cycle} 
\end{figure}

	The betweenness centrality \cite{freeman1978centrality} of a node $k$ is defined as the fraction of times that a node $i$, needs a node $k$ in order to reach a node $j$ via the shortest path. Specifically, if $g_{ij}$ is the number of shortest paths from $i$ to $j$, and $g_{ikj}$ is the number of $i$-$j$ shortest paths that use $k$, then the betweenness centrality of node $k$ is given by
	\begin{equation*}
	\sum\limits_i \sum\limits_j \frac{g_{ikj}}{g_{ij}}, \quad i \neq j \neq k,
\end{equation*}
	which essentially counts the number of shortest paths that pass through a node $k$ since we assume that $g_{ij} = 1$ in our road network because edges are weighted according to time. For a given road segment in the road segment graph, the betweenness would basically be the road segments share of all shortest-paths that utilize that the given road segment.  Intuitively, if we are given a road network containing two cities separated by a bridge, the bridge will likely have high betweenness centrality. It also seems like a good installation location because of the importance it plays in the network. Thus, for a small budget, we can expect the solution based on the betweenness centrality to give to be reasonable in such scenarios.  There is however an obvious downfall to this heuristic, consider a road network where the betweenness centrality of all the nodes are identical. For example, take a the cycle on $n$ nodes. Then using this heuristic would be equivalent to choosing installation locations at random. A cycle on $n$ vertices can represent a route taken by a bus, thus, a very practical example.  Figure \ref{cycle} shows an optimal solution from our model to minimize the number of infeasible routes with a budget of at most four units. 
The eigenvector centrality \cite{newman2008mathematics} of a network is also considered. As an extension of the degree centrality, a centrality measure based on the degree of the node, the concept behind the eigenvector centrality is that the importance of a node is increased if it connected to other important nodes. In terms of a road segment graph, this would translate into the importance of a road segment increasing if its adjacent road segments are themselves important. For example, if a road segment is adjacent to a bridge. One drawback of using this centrality measure is that degree of nodes in road segment graphs  is typically small across the graph. However, it stll helps to find regions of potentially heavy traffic.

	\section{Results and Discussion}
	In this section, the results of the proposed model for the WCL installation problem are discussed. In the following experiments Pyomo, a collection of Python packages \cite{hart2011pyomo,hart2012pyomo} was used to model the integer program. As a solver, CPLEX 12.7 \cite{ilog2014cplex} was used with all  results attained with an optimality gap of at most 10\%. Designing a fast customized solver is not the central goal of this paper. However, it is clear that introducing customized parallelization and using advanced solvers will make the proposed model solvable for the size of a large city in urban area.
	
	 The measurement of WCL installation effectiveness on a particular road segment depends on the $\soc$ function used. However, the $\soc$ function varies from EV to EV and is dependent on such factors as vehicle and battery type and size together with the effectiveness of the charging technology used. However, the purpose of this paper is to propose a model that is able to accommodate any $\soc$ function. 
\subsection{Small networks}     
\begin{figure}[ht]
\begin{center}
	 \subfigure[]{\includegraphics[width=0.33\linewidth]{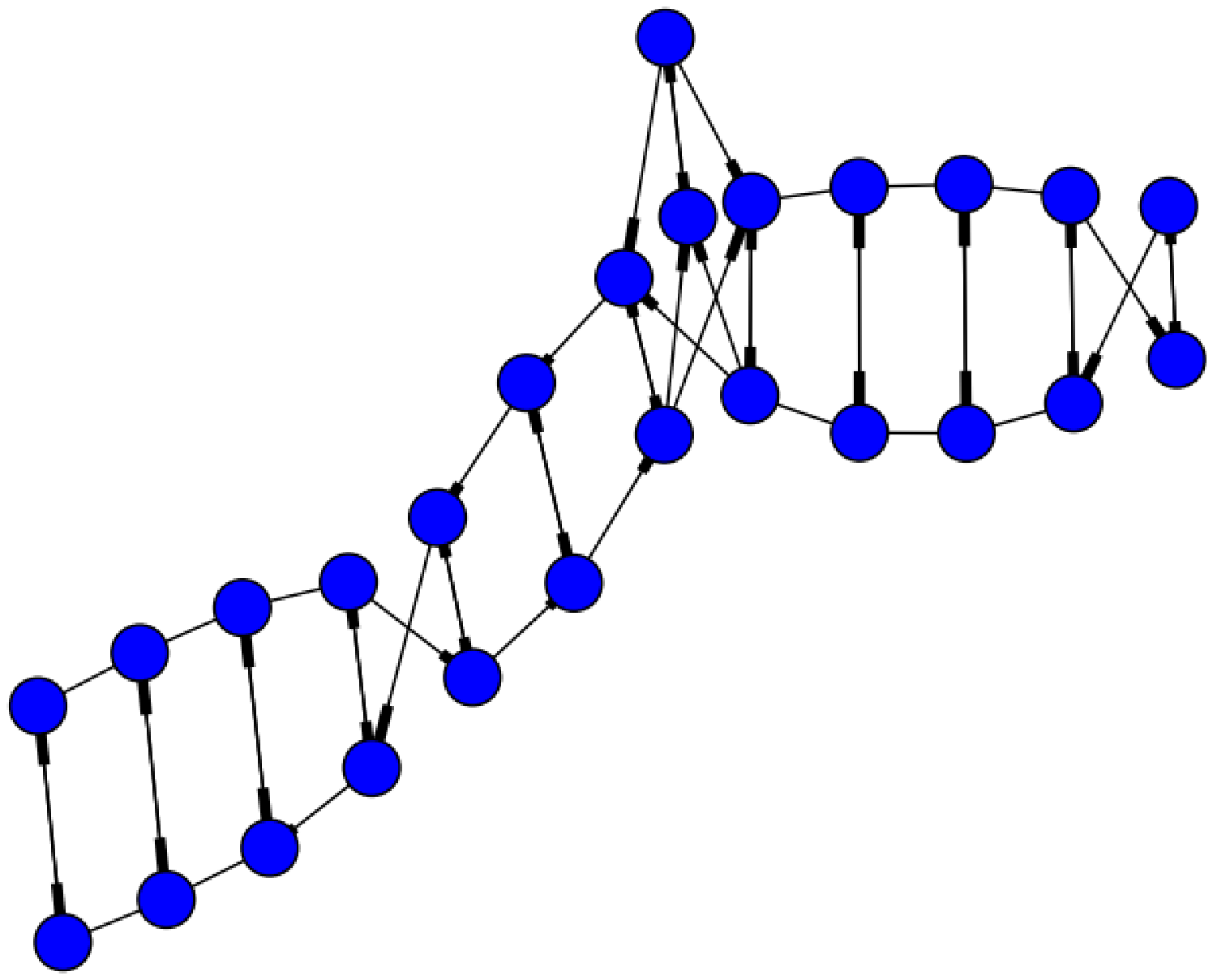}}
	\subfigure[]{ \includegraphics[width=0.33\linewidth]{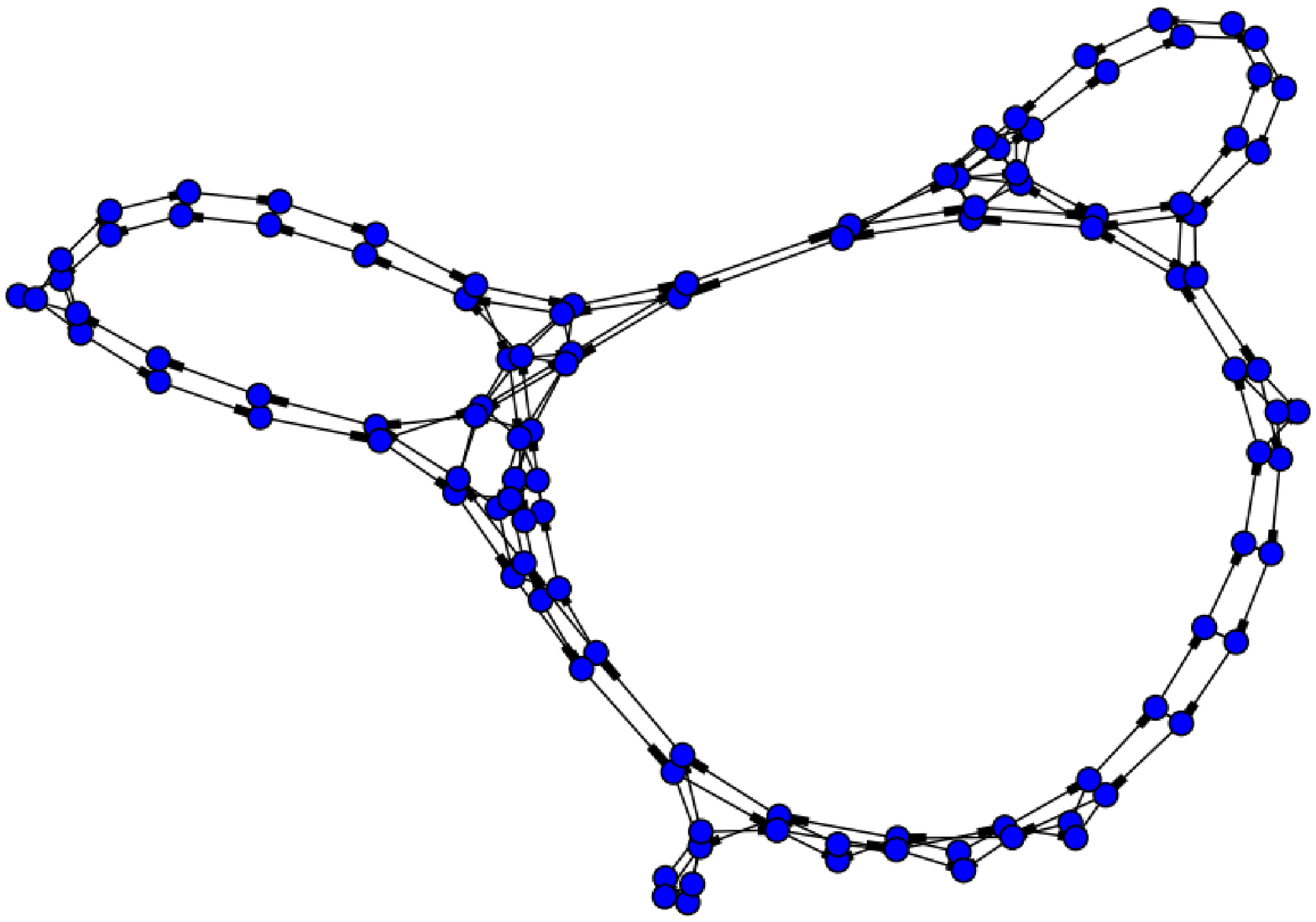}}	 
	 \caption[]{Directed toy graphs of 26 and 110 vertices used for problems 1 and 2, respectively. The bold end points on the edges of (a) represent edge directions. The graphs are subgraphs of the California road network taken from the dataset SNAP in \cite{leskovec2009community}
	 	 }	
	 \label{toygraph}
\end{center}

	\end{figure}

     In order to demonstrate the effectiveness of our model, we begin with presenting the results on two small toy graphs in which all road segments are identical. We incorporate a simplistic $\soc$ function, one in which the $\soc$ increases and decreases by one \soc discrete level if a wireless charging lane is installed or not installed, respectively. For the first toy graph (see Figure \ref{toygraph}(a)), the objective is to determine the minimum budget such that all routes are feasible. In this case, we assume that a fully charged battery has four different levels 
of charge 0, 1,2, and 3, where a fully charged battery contains three units. This would imply that the \soc-state graph would contain four levels. The parameter $\alpha$ is fixed to be 0. For the second toy graph (see Figure \ref{toygraph}(b)), the objective is to minimize the number of infeasible routes with a varying budget.  
 In this case, we take the number of discrete \soc levels in the \soc-state graph is taken to be five, with $\alpha=0$. 

\begin{figure}[ht]
		\begin{center}
			\subfigure[Optimal, $B=12$]{\includegraphics[width=0.325\linewidth]{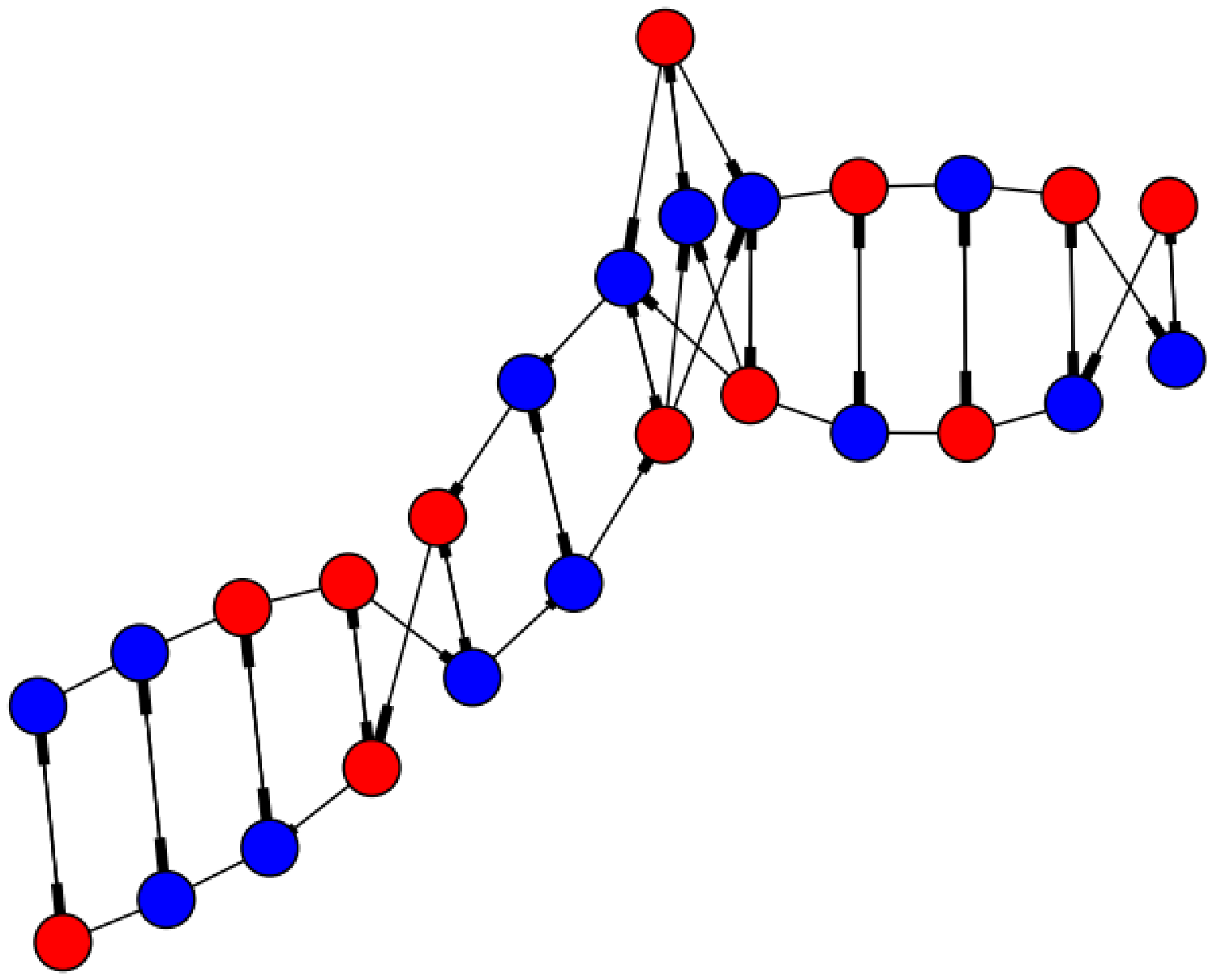}}
			\subfigure[Betweenness, $B=20$]{ \includegraphics[width=0.325\linewidth]{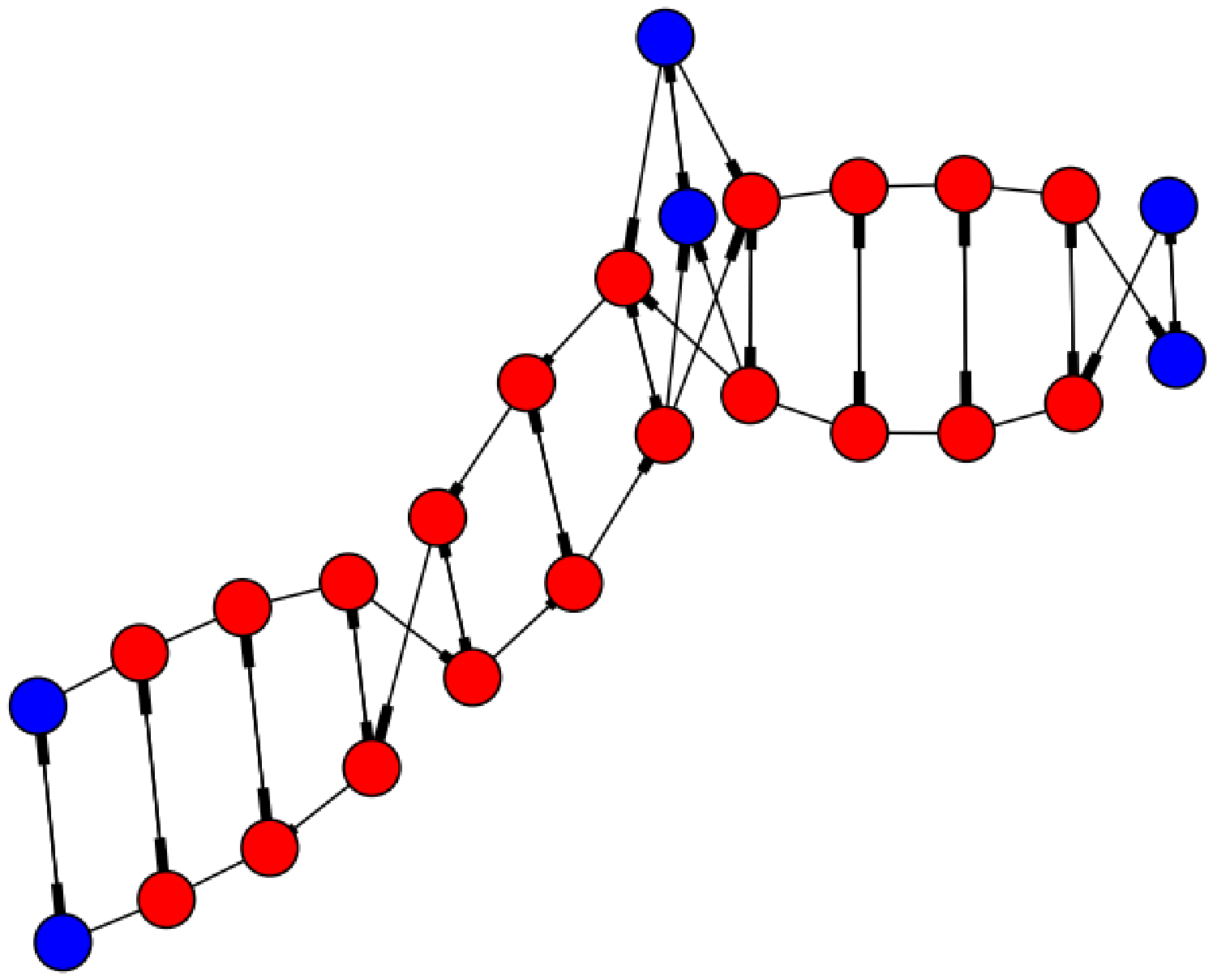}}
			\subfigure[Eigenvector, $B=23$]{ \includegraphics[width=0.325\linewidth]{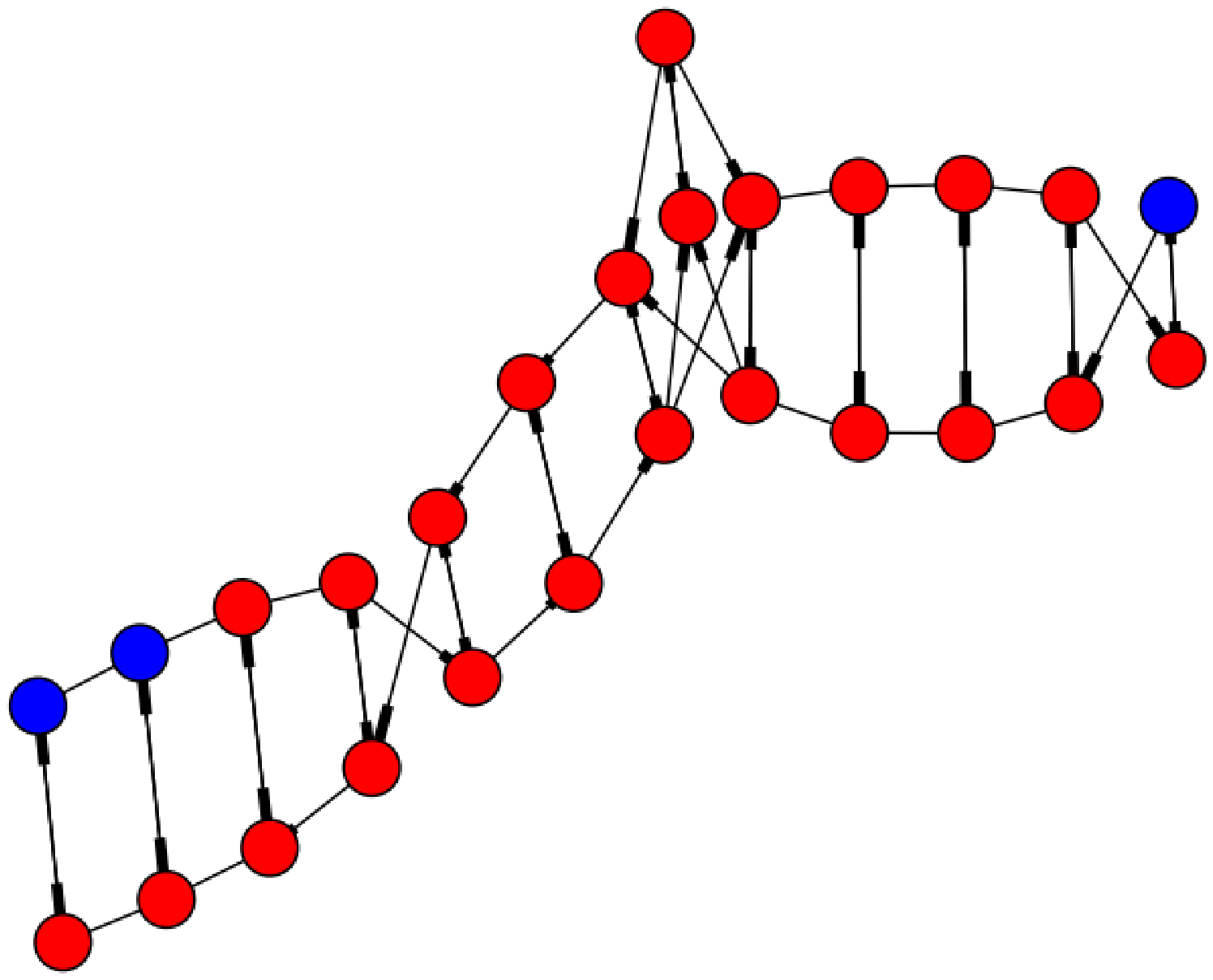}}
			\caption{Comparison of the different methods. The minimum number of WCL installation needed to eliminate all infeasible routes is $B$. The nodes colored red indicate location of WCL installation. In (a), we demonstrate the result given by our model requiring a budget of 12 WCL'sin order to have zero infeasible routes. In (b) and (c), we demonstrate  solutions from the betweenness and eigenvector heuristics that give budgets of 20 and 23 WCL's, respectively.} 
			\label{optimal_vs_heuristics}
		\end{center}
		
	\end{figure}

In the experiments with the graph shown in Figure \ref{toygraph}(a), we take \textbf{all} routes into consideration and compute an optimal solution which is compared with  the betweenness and eigenvector centralities. We rank the nodes based on their centrality indexes, and take the smallest number of top $k$ central nodes that ensure that all routes are feasible. The installation locations for each method are shown in Figure \ref{optimal_vs_heuristics} of which the solution to our model uses the smallest budget. We observe a significant difference in the required budget to ensure feasibility of the routes (see values $B$ in the figure).

\begin{figure}[ht]
    \centering
		\subfigure[$\beta = \frac{10}{110}$]{\includegraphics[width=0.4\linewidth]{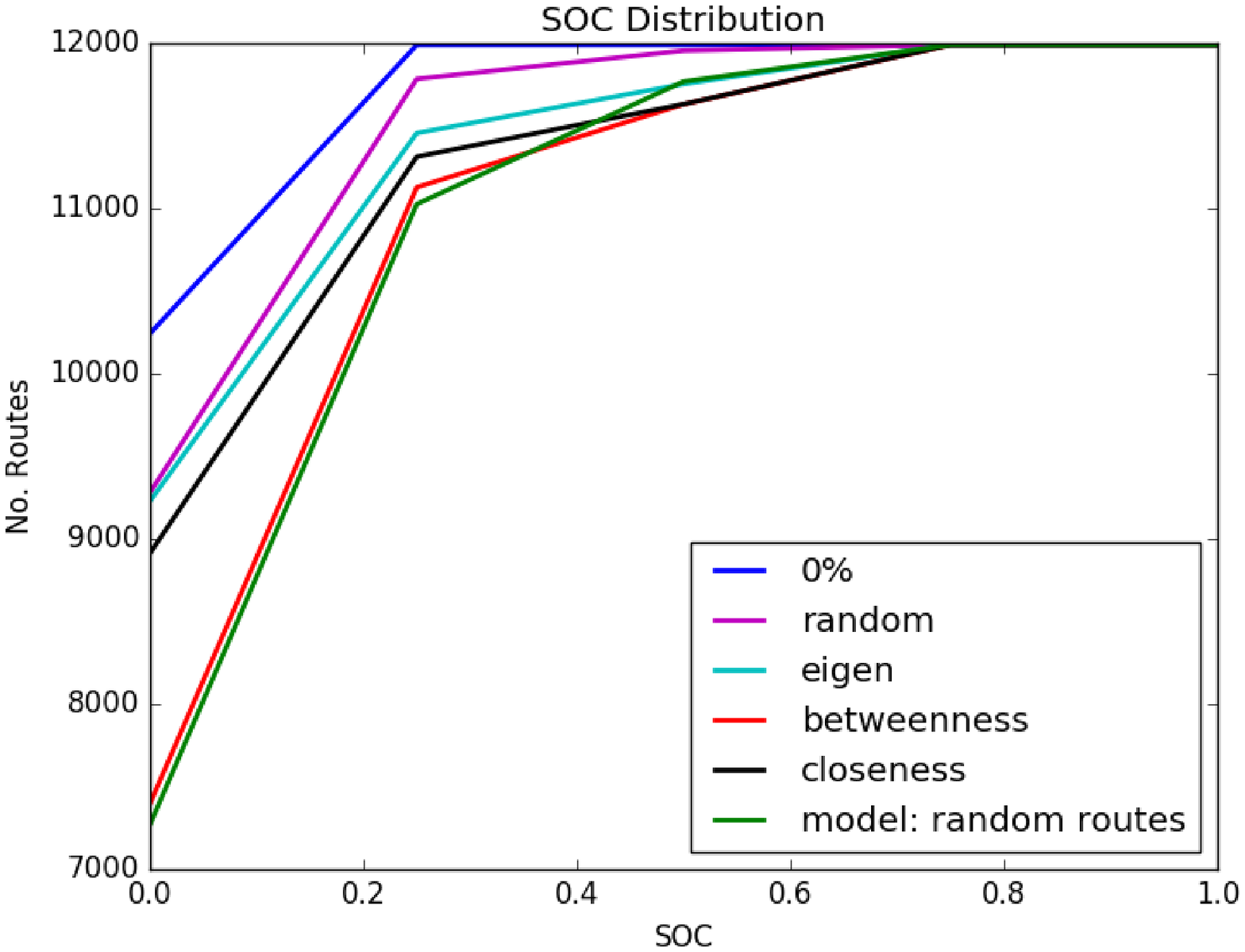}}
		\subfigure[$\beta = \frac{20}{110}$]{\includegraphics[width=0.4\linewidth]{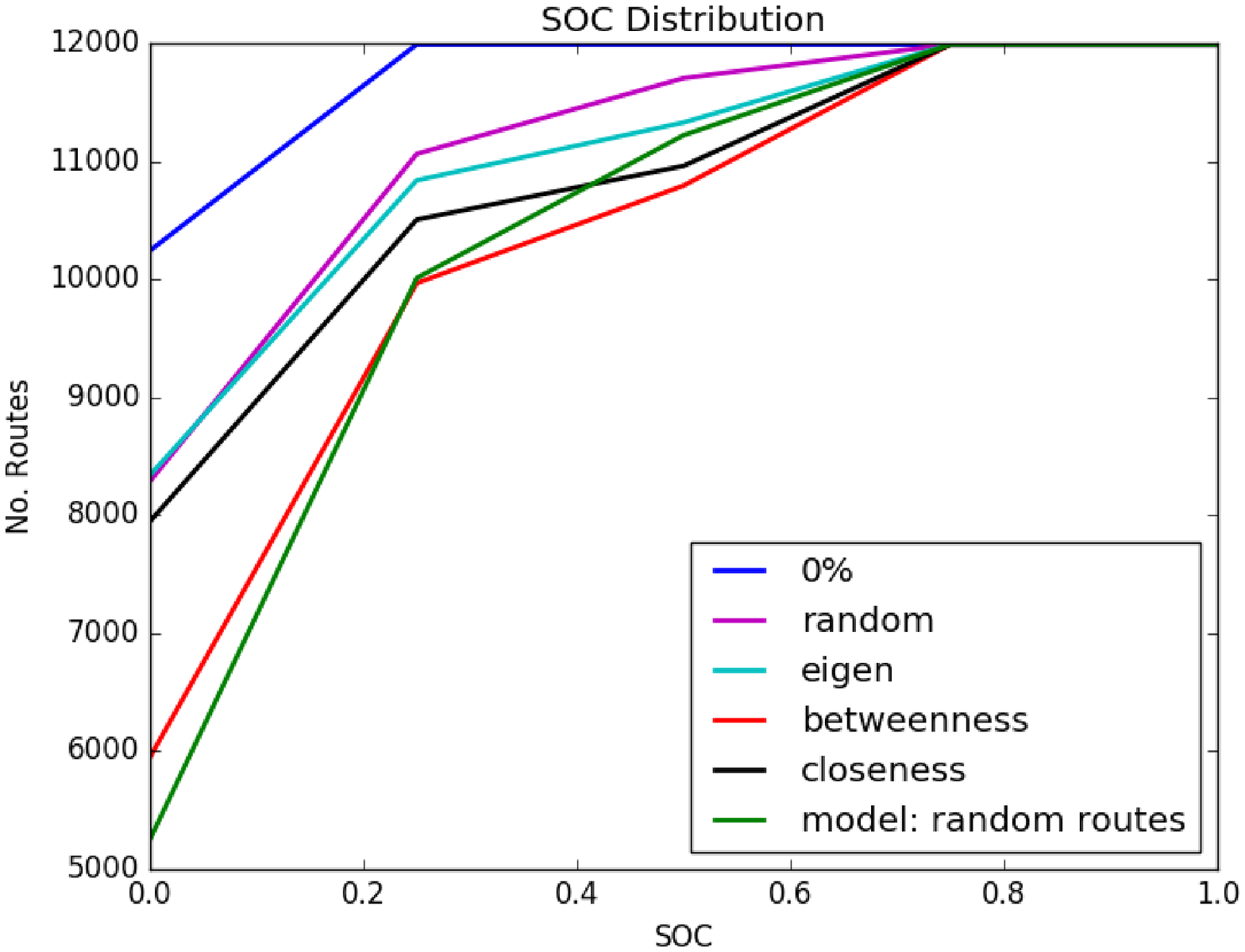}}
		\subfigure[$\beta = \frac{30}{110}$]{\includegraphics[width=0.4\linewidth]{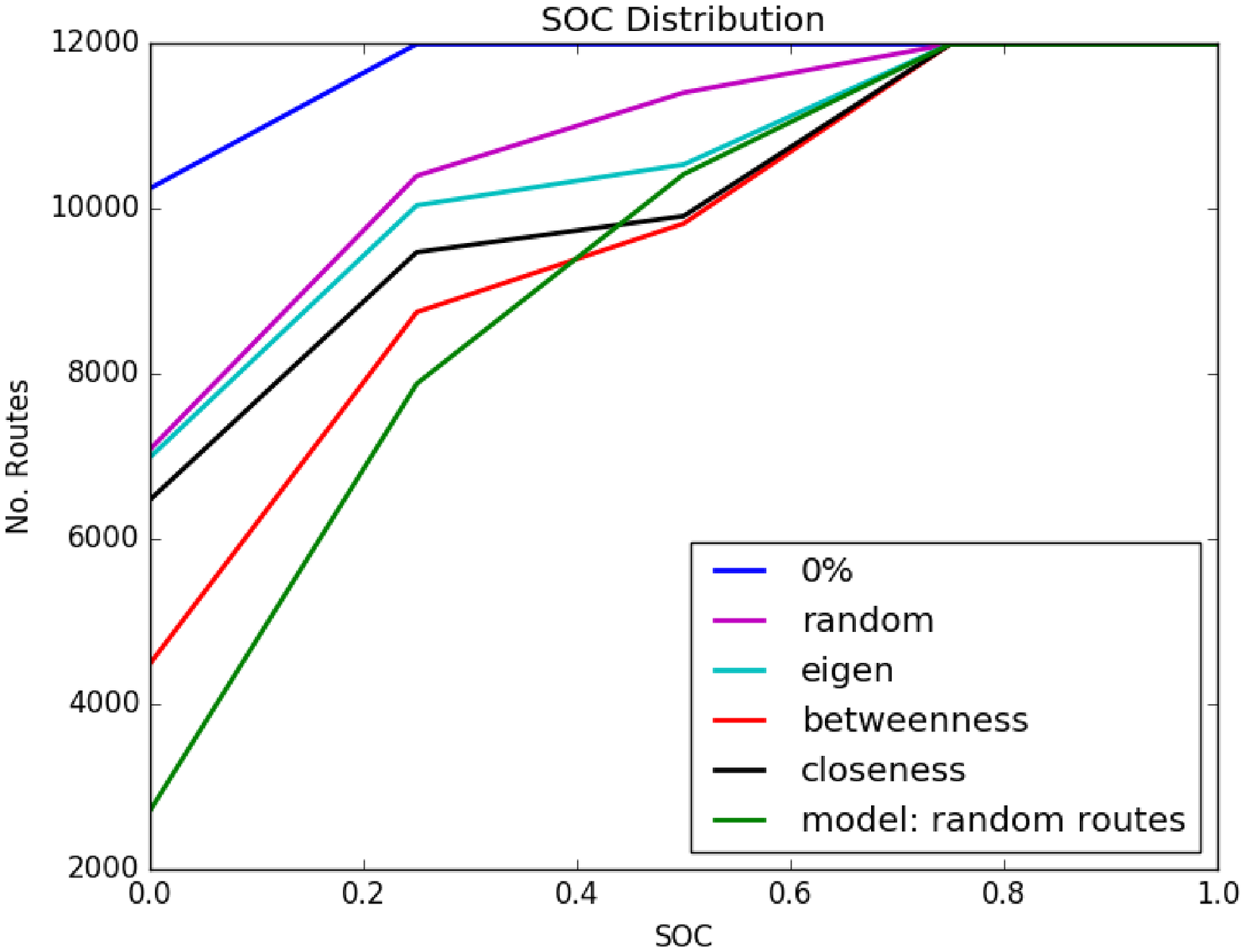}}
		\subfigure[$\beta = \frac{40}{110}$]{\includegraphics[width=0.4\linewidth]{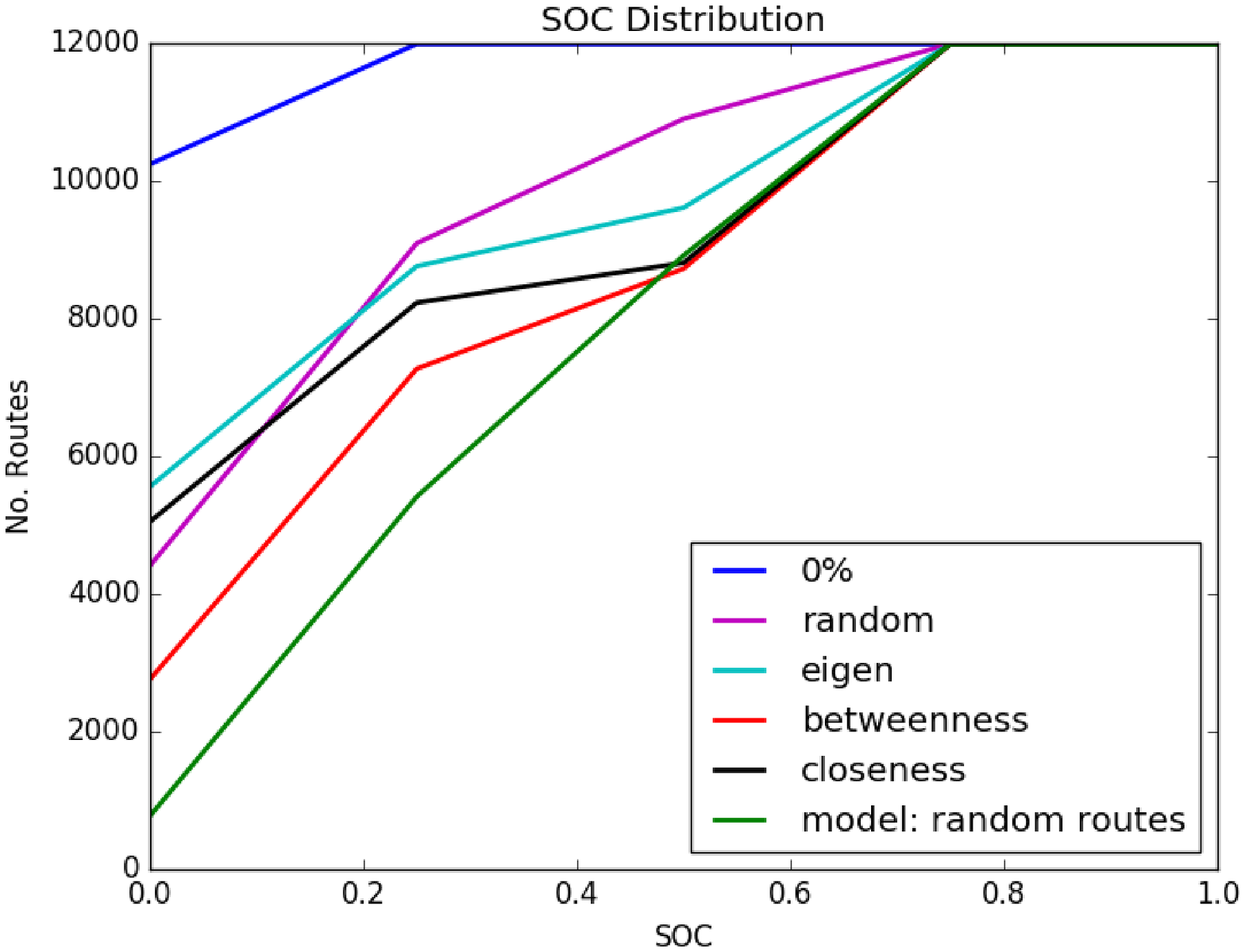}}
		\subfigure[$\beta = \frac{50}{110}$]{\includegraphics[width=0.4\linewidth]{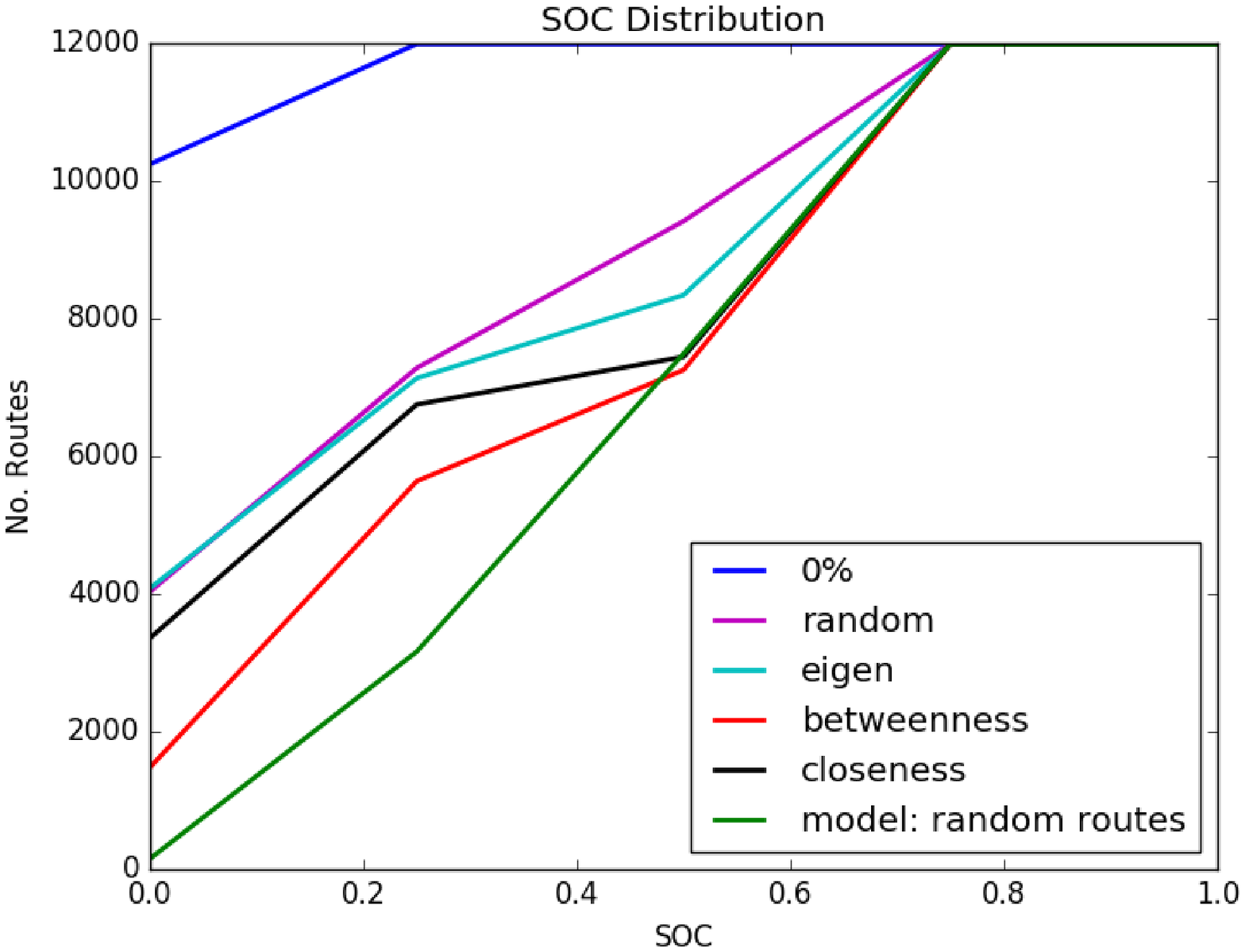}}
		\subfigure[$\beta = \frac{60}{110}$]{\includegraphics[width=0.4\linewidth]{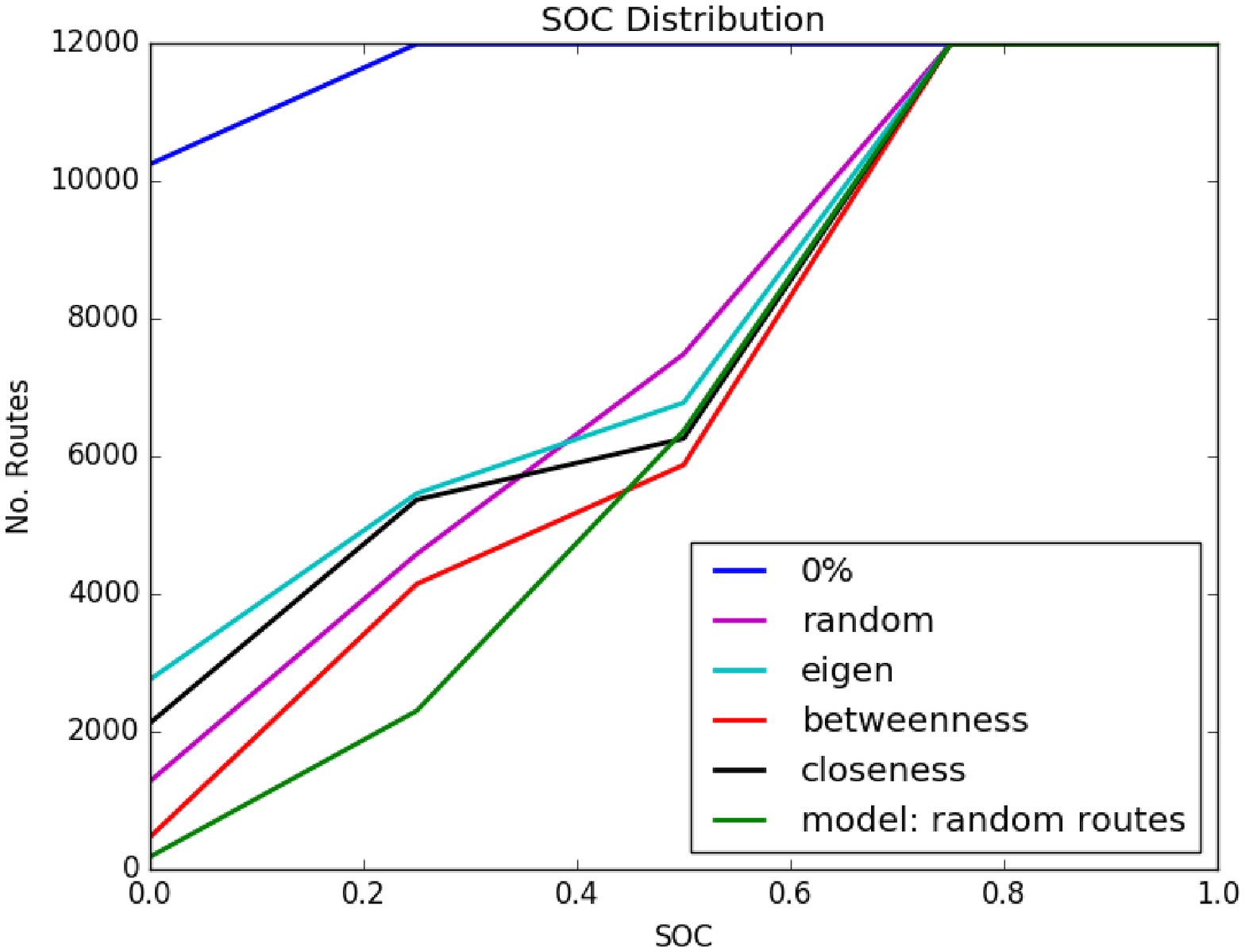}}
		\caption{Figures (a) to (f) are plots showing the number of routes ending with final $\soc$ below a given value via the different models. Legend ``model: random routes'' represents the solution from the proposed model when 100 routes were chosen uniformly at random with $\alpha=0$ and different budget scenarios. The solution is compared to the solutions from the different centrality measures, a random installation and one with no WCL installation. The $y$-intercept of the different lines shows the number of infeasible routes for the different methods. Our model gives a smaller number in all cases. The plots go further and show how a specific solution affects the $\soc$ of all routes. As the budget approaches 50\%, we demonstrate that our model gives a significant reduction to the number of infeasible routes while also improving the $\soc$ in general of the feasible routes }
		\label{toygraph_results}
	\end{figure}

For the graph shown in Figure \ref{toygraph}(b), we vary the available budget $\beta$. The results are shown in Figure \ref{toygraph_results}. The plots also indicate how the optimal solution affects the final $\soc$ of all other routes. The solutions to our model were based on 100 routes, with length at least 2, that were sampled uniformly without repetitions. We observe that our solution gives a very small number of infeasible routes for all budgets. We notice that for a smaller budget, taking a solution based on betweenness centrality gives a similar but slightly better solution than that produced by our model. However, this insignificant difference is eliminated as we increase the number of routes considered in our model. Note that if our budget was limited to one WCL, then the node chosen using the betweenness centrality would likely be a good solution because this would be the node that has the highest number of shortest paths traversed through it compared to other nodes. As we increase the budget, the quality of our solution is considerably better than the other techniques. For budgets close to 50\% in Figure \ref{toygraph_results} (d) and (e), 
our model gives a solution with approximately 90\% less infeasible routes compared to that of the betweenness centrality heuristic. This is in spite of only considering about 1\% of all routes as compared to betweenness centrality that takes all routes into account. 
\subsection{Experiments with Manhattan network}	
\begin{figure}[ht]
	\centering
\subfigure[Lower Manhattan]{\includegraphics[width=0.5\linewidth]{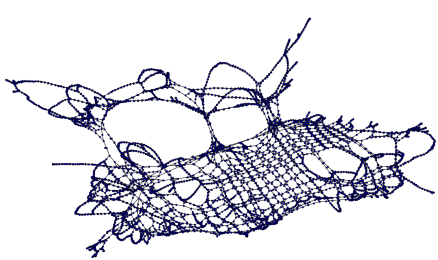}}
\subfigure[Manhattan Neighborhood]{\includegraphics[width=0.49\linewidth]{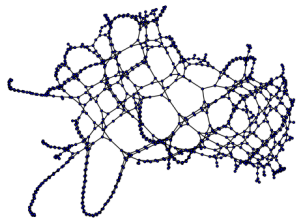}}
	\caption{Road segment graphs from real geospatial data: a node, drawn in blue, represents a road segment. Two road segments $u$ and $v$ are connected by a directed edge $(u,v)$ if and only if the end point of $u$ is that start point of $v$}	
	\label{manhattan}
\end{figure}

In the above example, the input to our model is a road segment graph with identical nodes, and a simple $\soc$ function. The proposed model was tested with real data and a realistic \soc function as defined in Section \ref{soc_section}. The data was extracted from lower Manhattan using OpenStreetMaps \cite{haklay2008openstreetmap}. The data was preprocessed by dividing each road into road segments. Each road in OpenStreetMaps is categorized into one of eight categories presented in Table \ref{road_category}, together with the corresponding speed limit for a rural or urban setting. For this work, roads from categories 1 to 5 were considered as potential candidates for installing wireless charging lanes due to their massive exploitation. Thus, any intersections that branch off to road categories 6 to 8 were ignored. The resulting road segment network contains 5792 nodes for lower Manhattan. Experiments on a neighborhood of 
lower Manhattan forming a graph of 914 nodes were also carried out.  The graphs are shown in Figure \ref{manhattan}. The authors would like to acknowledge the availability of other test networks such as the one maintained by \cite{bar2009transportation}.  
\begin{table}[htp]
	\centering
    	\caption{Road category with corresponding speed in Miles/Hr}
	\begin{tabular}{| l | l | l |l|}
		\hline
		Category & Road Type & Urban Speed & Rural Speed\\ \hline
		1&Motorway&60&70\\ \hline
		2&Trunk&45&55\\ \hline
		3&Primary&30&50\\ \hline
		4&Secondary&20&45\\ \hline
		5&Tertiary&15&35\\ \hline
		6&Residential/Unclassified&8&25\\ \hline
		7&Service&5&10\\ \hline
		8&Living street&5&10\\ \hline
	\end{tabular}
	\label{road_category}
\end{table}

Similar to experiments on the graph shown in Figure \ref{toygraph}(b) experiments are carried out on the Manhattan network using  200 routes. 
Routes that have a final $\soc$ less than the threshold $\alpha$ are sampled uniformly at random without repetitions. 
Due to relatively small driving radius within the Manhattan neighborhood graph shown in Figure \ref{manhattan} (b), 
the length of each road segment is increased by a constant factor in order to have a wider range of a final $\soc$ within each route. We take $\alpha = 0.8$ and 0.85 with a corresponding budget of $\beta=0.1$ and 0.2 respectively for the Manhattan neighborhood graph while $\alpha = 0.7$ and $\beta=0.1$ for the lower Manhattan graph. Results are compared with the heuristic of choosing installation locations based on their betweenness centrality. \emph{In our experiments, the betweenness centrality produces significantly better results than other heuristics, so it is used as our main comparison.} 

\begin{figure}[ht]
	\centering
	\subfigure[]{\includegraphics[width=0.329\linewidth]{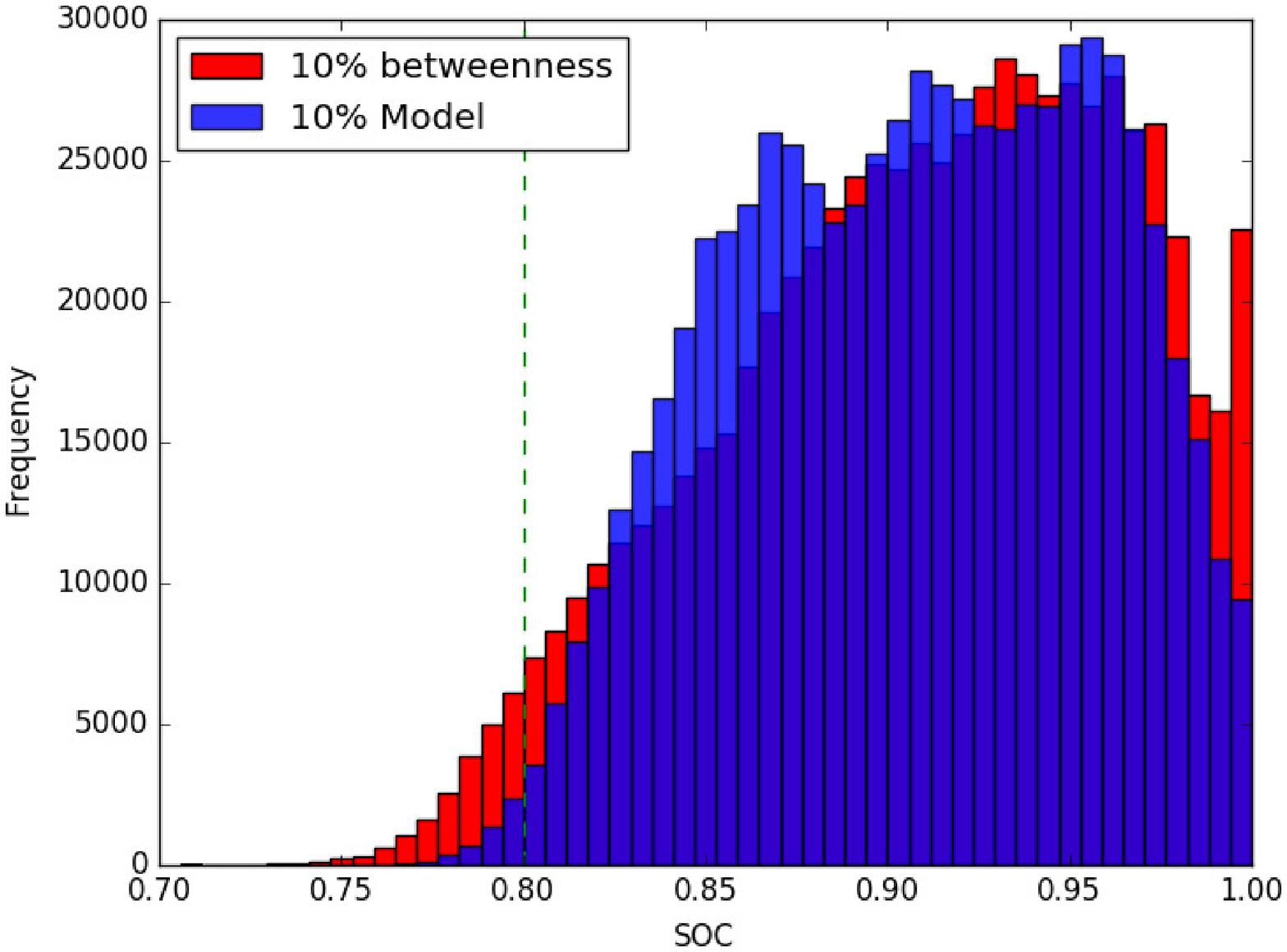}}
	\subfigure[]{\includegraphics[width=0.329\linewidth]{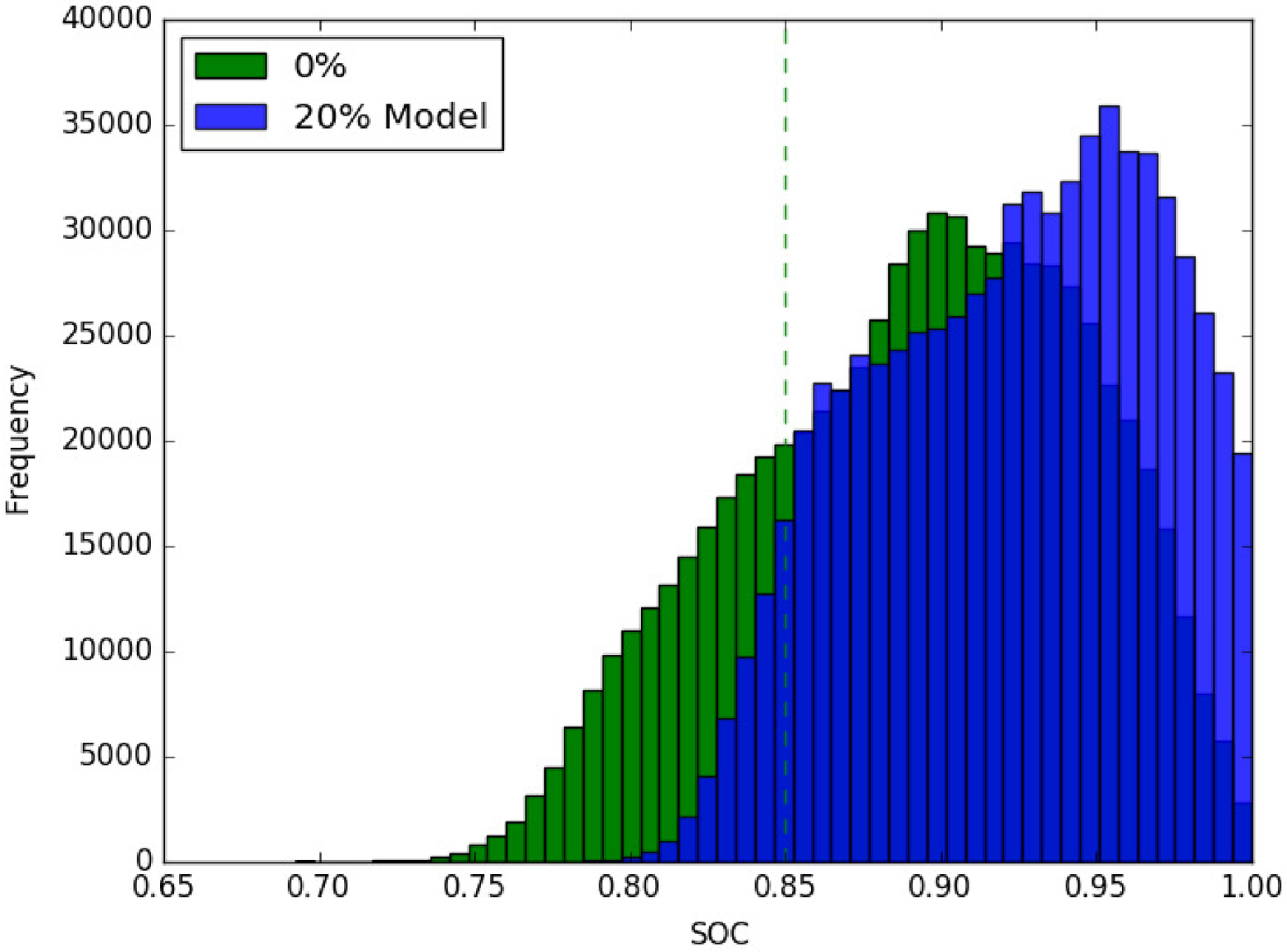}}
	\subfigure[]{\includegraphics[width=0.329\linewidth]{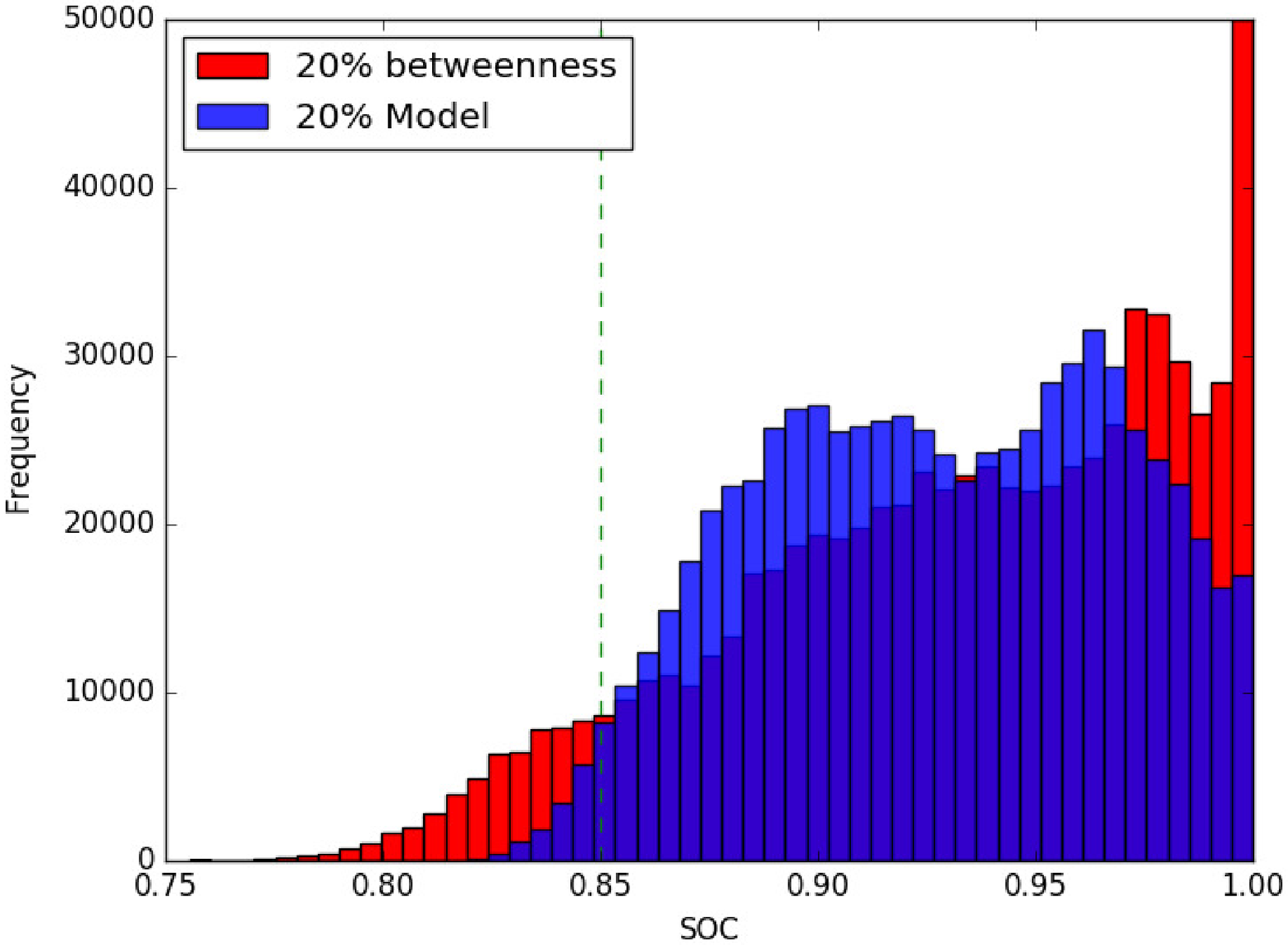}} 
\caption{
Histograms showing the number of infeasible routes for different values of $\alpha$ and $\beta$ for the Manhattan neighborhood graph. The vertical line indicates the value of $\alpha$. In (a) with a budget of 10\%, our model gives a solution with at least 50\% less infeasible routes compared to the betweenness heuristic. In (b), we demonstrate how the effects of a 20\% budget on the $\soc$ distribution within the network. In (c), our model gives a solution with at least 25\% less infeasible routes.}
	\label{neigh_results}
\end{figure}

For a threshold $\alpha=0.8$ in the Manhattan neighborhood graph, there are 42,001 infeasible routes with no WCL installation. With a budget $\beta=0.1$, the proposed model was able to reduce this number to 4,957. Using the heuristic based on betweenness centrality, the solution found contained 21,562 infeasible routes. For a budget of $\beta=0.2$ with threshold $\alpha=0.85$, there were 170,393 infeasible routes without a WCL installation, 57,564 using the betweenness centrality heuristic and only 14,993 using our model. Histograms that demonstrate the distributions of $\soc$ are shown in Figure \ref{neigh_results}. The green bars represent a $\soc$ distribution without any WCL installation. The red and blue bars represent  $\soc$ distributions after WCL installations based on the betweenness centrality heuristic and our proposed model, respectively.

\begin{figure}[ht]
\centering
\subfigure[]{\includegraphics[width=0.4\linewidth]{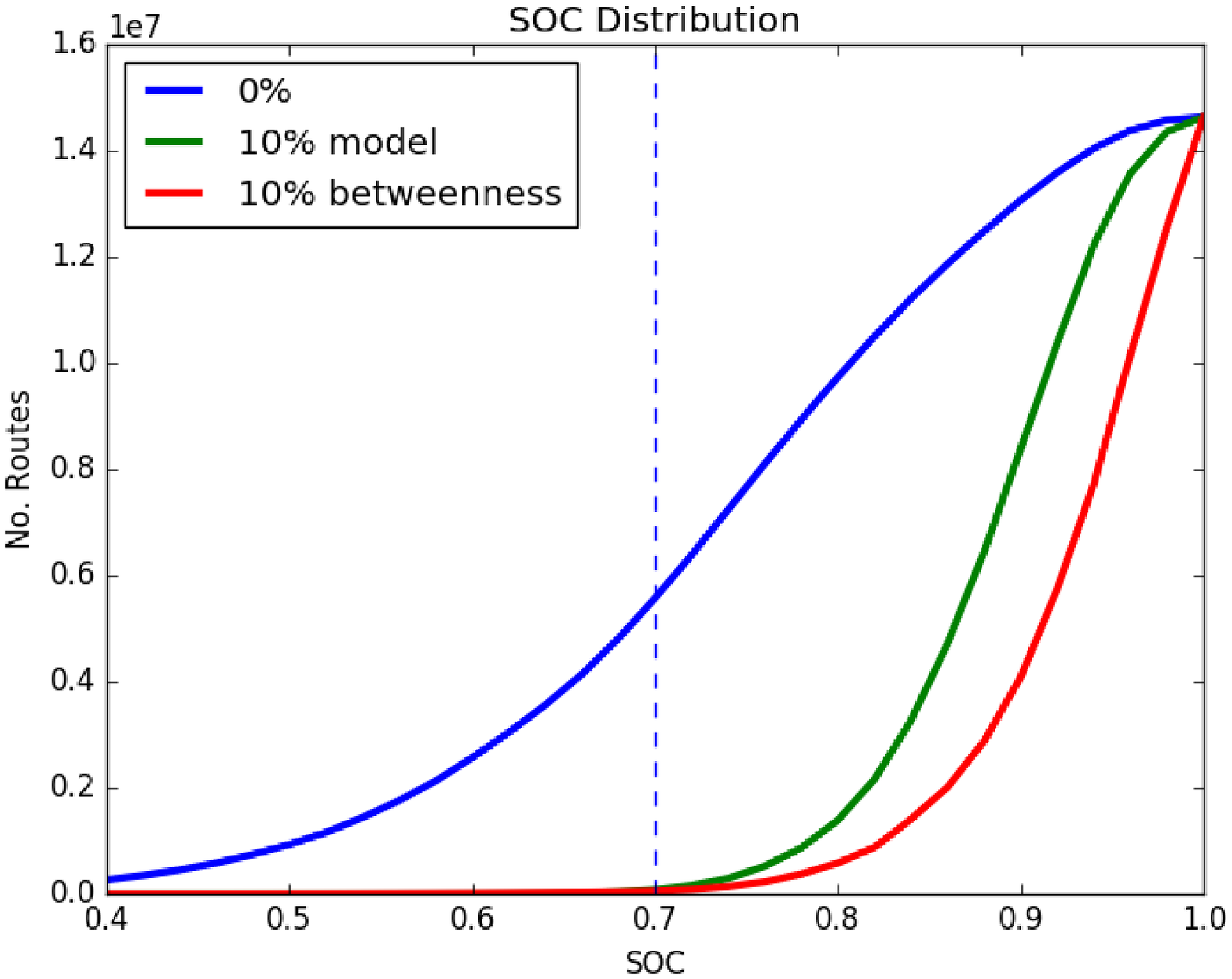}}
\subfigure[]{\includegraphics[width=0.4\linewidth]{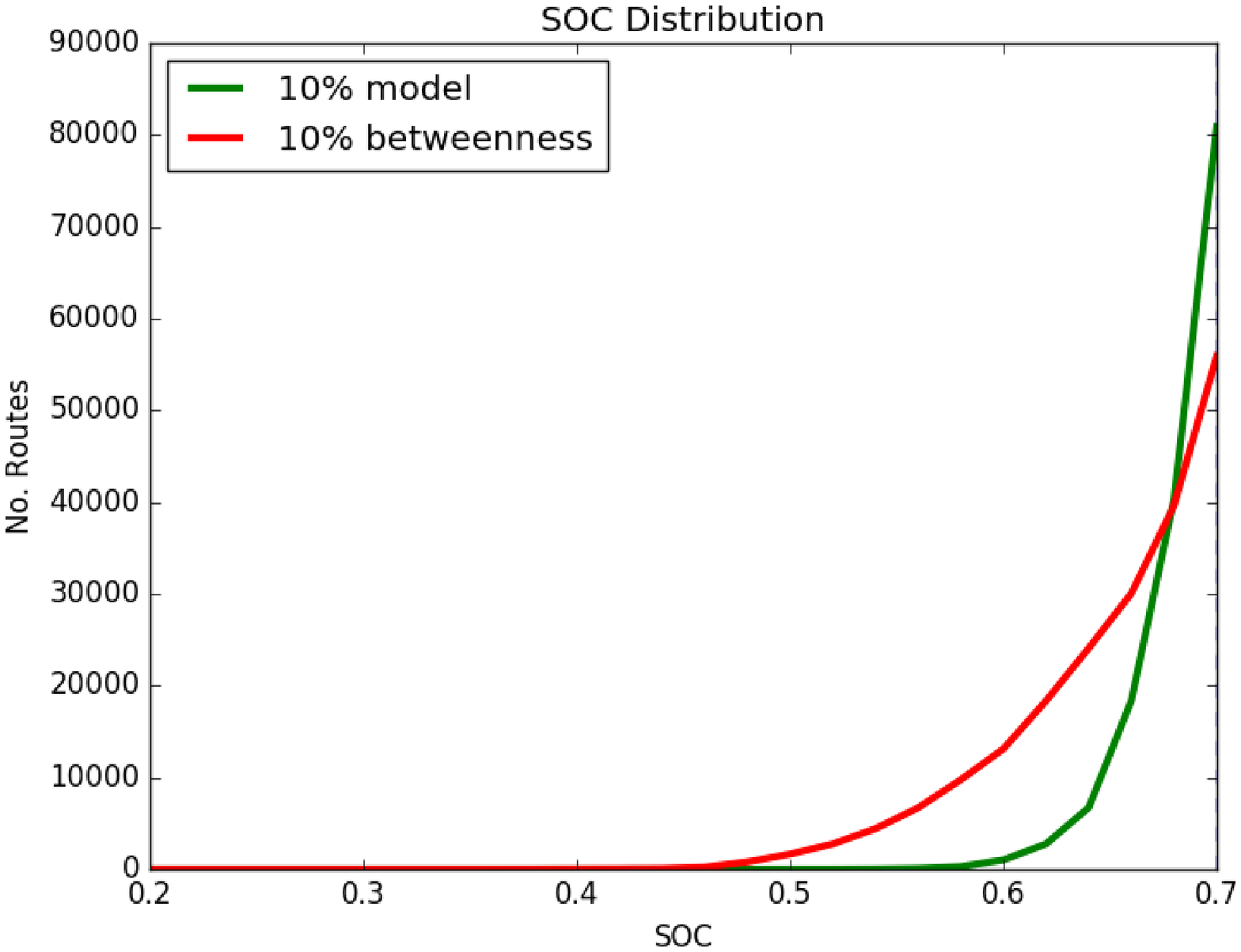}}
\caption{The number of routes ending with final $\soc$ below a given value in the lower Manhattan graph. The solution was obtained with $\alpha=0.7$ and $\beta=0.1$. The blue curve shows the $\soc$ distribution when no WCL are installed. Green and red curves show the $\soc$ distribution after an installation using the proposed model and the betweenness heuristic, respectively.  Plot (b) a gives closer look into (a) for the $\soc$ values below 0.7.}	
\label{manhattan_results}
\end{figure}

In Figure \ref{manhattan_results}, we demonstrate the results for the lower Manhattan graph, with $\alpha= 0.7$ and $\beta=0.1$. Due to the large number of routes, we sample about 16 million routes. From this sample, our model gives a solution with 10\% more infeasible routes compared to the heuristic based on betweenness centrality. Note that in this graph, there are about 13 million infeasible routes. From these routes, we randomly chose less than 1000 routes for our model without any sophisticated technique for choosing these routes, while the heuristic based on betweenness centrality takes all routes into account. The  plot in Figure \ref{manhattan_results}, shows the number of routes whose final $\soc$ falls below a given $\soc$ value. Similar to Figure \ref{toygraph_results} (a), the results demonstrate that for a relatively small budget, our model gives a similar result compared to the betweenness centrality heuristic.

\subsection{Experiments with Random Initial SOC} 
In the preceding experiments, EVs were assumed to start their journey fully charged. However, the assumption in our model was that the initial $\soc$ be any fixed value. Thus, as an alternative scenario, one can take the initial $\soc$ to follow a given distribution selected either by past empirical data or known geographic information about a specific area. For example, one could assume higher values in residential areas compared to non-residential areas. In this work, we carried out experiments where the initial $\soc$ was chosen uniformly at random in the interval $(a,1)$. The left endpoint of the interval was chosen such that the final $\soc$ associated to any route would be positive. In order to give preference to longer routes, we define $\Omega_l$ as the set of all routes with distance greater than $\tau + l\sigma$, where $\tau$ is the average distance of a route with standard deviation $\sigma$, for some real number $l$. We then study the average of the final $\soc$ of all routes in $\Omega_l$ which we denote as $\lambda_l$.

In the Manhattan neighborhood graph, Figure \ref{manhattan} (b), we chose 1200 routes as an input to the model. These routes were chosen uniformly at random from $\Omega_l$. In our solution analysis, we took $a=0.4$ and computed $\lambda_l$. 
Without any installation, we had $\lambda_l \approx 0.39$ for $l\geq 2$ while $0.5 \leq \lambda_l \leq 0.51$ based the betweenness heuristic with a installation budget of 20\% of the entire network. 
However, our model gives us $0.59 \leq \lambda_l \leq 0.71$ given the same 20\% installation budget. The value of $\lambda_l$ in this case significantly increases for an increase in $l$ or in the number of routes sampled from $\Omega_l$.

\subsection{Experiments with Variable Velocity}

\begin{figure}[ht]
\centering
\includegraphics[width=.5\linewidth]{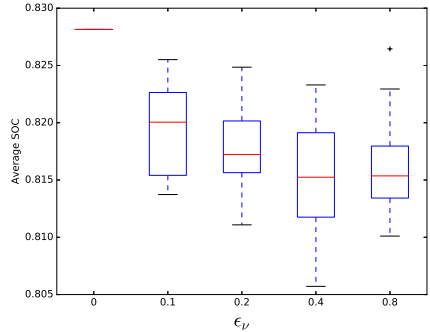}
\caption{Each boxplot depicts the average final \soc for 30 infeasible routes under 50 randomly generated velocity modification scenarios for each $\epsilon_{\nu}$. 
}
\label{velocityexp}
\end{figure}

The velocity of an EV along a road segment is assumed to be constant. We carry out experiments to determine how the solution of the proposed model may be affected in the case of a variation of the velocity of a road segment. The experiments consists of modifying the velocity along each road segment by at most a fixed percentage represented as $\epsilon_{\nu}\cdot100$ percent, where $\epsilon_{\nu} \in (-1,1)$.
The solutions derived after each road segment's velocity is modified is compared to the solution using the original velocity, in other words, $\epsilon_{\nu} =0$. 
For this experiment, we pick a set of 30 infeasible routes, with $\alpha = 0.8, \beta = 0.1$, chosen such that the distance of each route is at least $\tau   + 2\sigma$, where $\tau$ and $\sigma$ is the average distance of a route, and standard deviation respectively. For each $\epsilon_v$ we run 50 experiments computing the average final \soc of each run.  Figure \ref{velocityexp} summarizes the results which show that the quality of the solutions gradually decreases as more uncertainty is added into the velocity of each road segment. 

\subsection{Using Betweenness Centrality as Initial Solution}
Various studies \cite{he2013integrated,riemann2015optimal,chen2016optimal} have formulated different mathematical models, under varying assumptions, for the optimal placement of WCL's. In their work, due to the large number of variables generated, small size road networks  (usually under 80 road segments and 30 intersections) are considered. The network topology is implicitly taken into consideration in their mathematical model. In this work, the results presented in Figure \ref{toygraph_results} show that for a very small budget, heuristics from network analysis, such as ones based on betweenness centrality can give a reasonable solution in relation to their computational time. However, as expected, these heuristics are not  optimal. 

\begin{figure}[ht]
\centering
\includegraphics[width=.5\linewidth]{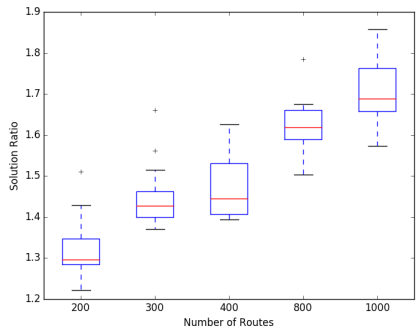}
\caption{Using Betweenness Centrality Heuristic to find an Initial Solution}
\label{btn_exp}
\end{figure}

In order to cut the computational time for solving problems on large-scale road networks, we carry out experiments where solutions from the heuristic based on betweenness centrality is used as an initial feasible solution to the problem. The experiment consists of selecting $k$ routes at random as input to the model.  For a fixed amount of time, we compare the objective values $obj_1$ and  $obj_2$ where $obj_1$ and $obj_2$ are the objective values as a result of the two starting conditions; 1) an initial solution is provided and 2) an initial solution is not provided, respectively. We run the experiment 30 times for each $k$, and vary $k$ between 200 and 1000 routes. In this experiment, we use the Manhattan neighborhood graph which consists of 914 road segments. Figure \ref{btn_exp} summarizes the results. The $y$-axis gives the ratio of the objective values $obj_1/obj_2$ measured after running CPLEX for at most 1 hour. The results show that using heuristics from network science can significantly lead to better results if the CPLEX is run for a fixed amount of time. For example, when 1000 routes were given as input to the model, and an initial solution based on betweenness centrality was provided, CPLEX, solving the maximization problem, recorded an objective value 70\% higher than one when an initial solution was not provided for over half of the test cases. This is a significant improvement for limited budget real settings.
\section{Conclusion and Future Work}
In this work, we have presented an integer programming formulation for modeling the WCL installation problem. With a  modification to the WCL installation optimization model, we present a formulation that can be used to answer two types of questions. First, determining a minimum budget to reduce the number of infeasible routes to zero, thus, assuring EV drivers of arrival at their destinations with a battery charge above a certain threshold. Second, for a fixed budget, minimizing the number of infeasible routes and thus reducing drivers' range anxiety.
 
    Our experiments have shown that our model gives a high quality solution that typically improves various centrality based heuristics. The best reasonable candidate (among many heuristics we tested) that sometimes not significantly outperforms our model is the betweenness centrality.
In our experiments, the routes were chosen randomly based on whether their final $\soc$ is below $\alpha$ or not. 
 We notice that a smarter way of choosing the routes leads to a better solution, for example choosing the longest routes generally provided better solutions. In future research, a careful study on the choice of routes to include in the model will give more insight into the problem. EV drivers may choose to slow down when using a charging lane to get more energy from the system or at intersections via opportunistic charging \cite{khan2017utility}. As a future direction, a probabilistic measure representing the chance of slowing down on a given route (for example, depending on the initial \soc) will be included in the model. For a more comprehensive study of this model and its desired modifications, we suggest to evaluate their performance using a large number of artificially generated networks using \cite{gutfraind2015multiscale,staudt2016generating}.

\section*{Acknowledgments} This research is supported by the National Science Foundation under Award \#1647361. Any opinions, findings, conclusions or recommendations expressed in this material are those of the authors and do not necessarily reflect the views of the  of the National Science Foundation. The authors are very grateful to the anonymous referees for their insightful comments. Clemson University is acknowledged for generous allotment of compute time on Palmetto cluster.

\bibliographystyle{AIMS}
\bibliography{mybib}
\medskip
\medskip

\end{document}